\documentclass[graybox]{svmult}

\usepackage{amsmath}
\usepackage{mathptmx}       
\usepackage{helvet}         
\usepackage{courier}        
\usepackage{type1cm}        
\usepackage{makeidx}         
\usepackage{graphicx}        
\usepackage{multicol}        
\usepackage[bottom]{footmisc}

\makeindex             
\usepackage{amssymb,epsfig}
\usepackage{color}
\usepackage{epsfig}
\usepackage{algorithm, algpseudocode}
\usepackage{subfigure}
\usepackage{multirow}
\usepackage{pgfplots}
\definecolor{dgreen}{RGB}{0,130,150}
\begin{document}
\title*{Efficient Resolution of Anisotropic Structures}
\titlerunning{Anisotropic Structures}
\author{W. Dahmen, C. Huang, G. Kutyniok, W.-Q Lim, C. Schwab, G. Welper}
\authorrunning{W. Dahmen, C. Huang, G. Kutyniok, W.-Q Lim, C. Schwab, G. Welper}
\institute{Wolfgang Dahmen \at IGPM, RWTH Aachen
\email{dahmen@igpm.rwth-aachen.de} \and Chunyan Huang \at School of
Applied Mathematics, Central University of Finance and Economics,
\email{hcy@cufe.edu.cn} \and Gitta Kutyniok \at Department of
Mathematics, {Technische Universit\"at Berlin} \email{kutyniok@math.tu-berlin.de} \and
Wang-Q Lim \at Department of Mathematics, {Technische Universit\"at Berlin}
\email{{lim@math.tu-berlin.de}} \and Christoph Schwab \at ETH Z\"urich
\email{christoph.schwab@sam.math.ethz} \and Gerrit Welper \at
Department of Mathematics, Texas A \& M University
\email{welper@math.tamu.edu} }
%
%
\maketitle

\vspace*{-3.5cm}
\begin{abstract}\newline
We highlight some recent new developments concerning the
{\em sparse representation} of possibly {\em high-dimensional}
functions exhibiting strong {\em anisotropic} features and low regularity
in isotropic Sobolev or Besov scales.
Specifically, we focus on the solution of {\em transport equations}
which exhibit {\em propagation of singularities}
where, additionally,
high-dimensionality enters when the convection field,
and hence the   solutions,
depend on parameters varying over some compact set.
Important constituents of our approach are directionally adaptive
discretization concepts motivated by compactly supported shearlet systems,
and well-conditioned stable variational formulations
that support  trial spaces with anisotropic refinements with arbitrary directionalities.
We prove that they provide tight error-residual relations which are used
to contrive rigorously founded adaptive refinement schemes which converge in $L_2$.
Moreover, in the context of parameter dependent problems
we discuss two approaches serving different purposes and working under different
regularity assumptions.
For ``frequent query problems'', making essential use
of the novel well-conditioned variational formulations,
a new {\em Reduced Basis Method} is outlined which exhibits a
certain {\em rate-optimal} performance
for indefinite, unsymmetric or singularly perturbed problems.
For the radiative transfer problem with scattering a {\em sparse tensor} method
is presented which mitigates or even overcomes the curse of dimensionality
under suitable (so far still isotropic) regularity assumptions.
Numerical examples for both methods illustrate the theoretical findings.
\end{abstract}

\vspace*{-0.1cm}


%
\vspace*{-4mm}
\setcounter{section}{0}
\section{Introduction}\label{dahmen_huang_kutyniok_lim_schwab:sect1}
%
The more complex a data site or mathematical model is the more adapted a corresponding mathematical representation needs to be in order to
capture its information content
at acceptable cost in terms of storage and computational complexity. In principle, this is true for mathematical objects described explicitly
by large sets of possibly noisy or corrupted data but also for those given only
implicitly as the solution of an operator equation.
The latter scenario is perhaps even more challenging because direct observations are not possible.
By ``adapted representation'' we mean a representation of
the unknown function that exploits possibly global features of this function
so as to require, for a prescribed target accuracy,
 only relatively few parameters to determine a corresponding approximation.
Such global features could take a variety of forms such as (i) a
high degree of regularity
except at isolated singularities {\em located on lower dimensional manifolds}, or
(ii) a particular {\em sparsity} possibly  with respect
to a {\em dictionary} which may even depend on the problem at hand.
In fact, corresponding scenarios are not strictly disjoint.
In either case reconstruction or approximation methods are necessarily nonlinear.
For instance, as for (i), 1D {\em best $N$-term} wavelet approximations offer a powerful method based on selecting only possible few coefficients in an exact representation with respect
to a given {\em universal background} dictionary, e.g. a wavelet basis. When dealing with more than one spatial variable
the situation quickly becomes more complicated and for  spatial dimensions much larger than three, classical numerical tools designed for the low dimensional regime become practically
useless.
This is commonly referred to as {\em curse of dimensionality}.
Unfortunately, there seems to be no universal strategy of
dealing with the curse of dimensionality, i.e.,
that works in all possible cases.

One global structural feature which is encountered in many multivariate
scenarios is {\em anisotropy}:
images, as fuctions of two variables, exhibit edges and discontinuities along curves.
Higher dimensional biological images have sharp interfaces separating
more homgeneous regions.
Likewise highly anisotropic phenomena such as shear- or boundary layers are
encountered in solutions to transport dominated initial-boundary value problems.

One major focus of this project has been to efficiently recover and
economically encode {\em anisotropic} structures
represented by explicitly given data or determined as solutions of operator equations which are
prone to give rise to such structures. Regarding this latter case, which we will focus on in this article,
{\em parametric transport problems} (as well as close relatives)
have served as guiding model problems for the following reasons: (i) their solutions could exhibit shear or boundary layers and hence
discontinuities across lower dimensional manifolds calling for suitable {\em anisotropic discretizations};
(ii) how to contrive suitable  {\em variational formulations}, which in particular accommodate such anisotropic discretizations is much less  clear than in the elliptic case;
 (iii) {\em parametric versions} give rise to {\em high-dimensional} problems.

 Concerning (i),  {\em directional representation systems} like {\em curvelets} and {\em shearlets}
 {outperform} classical {\em isotropic} wavelet bases when approximating so called ``cartoon images'',
 see  {\cite{dahmen_huang_kutyniok_lim_schwab:Do01}} and {\cite{dahmen_huang_kutyniok_lim_schwab:CaDo02,dahmen_huang_kutyniok_lim_schwab:KuLeLi12,dahmen_huang_kutyniok_lim_schwab:KuLi11,dahmen_huang_kutyniok_lim_schwab:Li10}}.
For recent applications  to imaging data, in particular, {\em inpainting} as well as in combination with
{\em geometric separation concepts} the reader is referred to  \cite{dahmen_huang_kutyniok_lim_schwab:DoKu13,dahmen_huang_kutyniok_lim_schwab:KiKuZh13}.
In the present context of solving operator equations
we outline in Section \ref{dahmen_huang_kutyniok_lim_schwab:sect2}
 trial spaces which accommodate directional adaptivity.
They are motivated by recent constructions
of  compactly supported {\em piecewise polynomial shearlet systems}
(see e.g. \cite{dahmen_huang_kutyniok_lim_schwab:KiKuLi12}) because they are close to classical multiresolution structures
and similar in nature to classical discretization systems.
Since cartoons exhibit structural similarities with the solution to transport problems
we state best $N$-term error bounds
for cartoon functions that will later serve as benchmarks for an adaptive solver.
For related anisotropic simplicial discretizations
and their analysis see e.g. \cite{dahmen_huang_kutyniok_lim_schwab:ChSuXu07,dahmen_huang_kutyniok_lim_schwab:CoDyHeMi12,dahmen_huang_kutyniok_lim_schwab:CoMi12}.

As for (ii), our approach differs from previous works on anisotropic discretizations derived from ``curvature information'' on the current approximation and hence  not based on a rigorous error control (see e.g. \cite{dahmen_huang_kutyniok_lim_schwab:Do98}
and the references therein),
in that we derive first  in Section \ref{dahmen_huang_kutyniok_lim_schwab:stabvar} {\em well conditioned variational formulations} for general unsymmetric or indefinite and singularly
perturbed problems, see \cite{dahmen_huang_kutyniok_lim_schwab:CoDaWe12,dahmen_huang_kutyniok_lim_schwab:DaHuScWe12} for details on {\em convection-diffusion} and {\em transport problems}.
The underlying basic principles are of independent interest by themselves and seem to have appeared first in \cite{dahmen_huang_kutyniok_lim_schwab:BaMo84}.
They are also closely related to ongoing developments running under
the flag of {\em Discontinuous Petrov Galerkin (DPG) Methods},
see e.g. \cite{dahmen_huang_kutyniok_lim_schwab:DeGo10,dahmen_huang_kutyniok_lim_schwab:DeGo11}.
The approach is motivated by two crucial corner stones.
On the one hand, one can essentially {\em choose} the norm for the (infinite dimensional) trial space $X$ by which one would like to measure accuracy
while adapting the norm for the (infinite dimensional) test space $Y$ so as to ensure that (ideally) the operator induced by this variational
formulation is even an isometry from $X$ to $Y'$ (the normed dual of $Y$).
Numerical feasibility of (nearly optimal)
Petrov Galerkin discretizations  based on such formulations,
even beyond a DPG framework, hinges on an appropriate {\em saddle point formulation}
which turns out to be actually  crucial in
connection with {\em model reduction} \cite{dahmen_huang_kutyniok_lim_schwab:DePeWo13}.
On the one hand, this allows one to accommodate, for
instance, $L_2$-frames.
On the other hand, the resulting tight error-residual relation
is the basis of computable a-posteriori error estimators \cite{dahmen_huang_kutyniok_lim_schwab:CoDaWe12,dahmen_huang_kutyniok_lim_schwab:DaHuScWe12} and,
ultimately, to rigorously founded adaptive anisotropic refinement strategies.

These variational formulations apply in much more generality but in order to address issue (iii) we exemplify them
for the simple linear transport equation (stationary or instationary) whose {\em parametric} version leads to high-dimensional problems
and forms a core constituent of kinetic models such as  {\em radiative transport}.
There the transport direction - the parameter - varies over a unit sphere so that solutions
are functions of the spatial variables (and, possibly, of time) and of the transport direction.

We briefly highlight two ways of treating such parametric problems under slightly different objectives.
Both strategies aim at approximating the solution
$u(x,\vec{s})$, $x\in \Omega \subset \mathbb{R}^d$, $\vec{s}\in S^{d-1}$, in the form\vspace*{-2mm}
\begin{equation}
\label{dahmen_huang_kutyniok_lim_schwab:separation}
u(x,\vec{s}) \approx \sum_{j=1}^n c_j(\vec{s}) u_j(x).
\end{equation}
%
In Section \ref{dahmen_huang_kutyniok_lim_schwab:RBM} the $u_j$ are constructed {\em offline} in a greedy manner from {\em snapshots} of the solution manifold, thus forming
a {\em solution dependent} dictionary. According to the paradigm of the {\em Reduced Basis Method} (RBM) the parameter dependent coefficients
$c_j(\vec{s})$ are not given explicitly but can be efficiently computed in an {\em online} fashion, e.g. in the context of design or
(online) optimization. This approach works the better the smoother
the dependence of the solution on the parameters is so that the Kolmogorov $n$-widths decay rapidly with increasing $n$.
Making essential use of the well conditioned variational formulations from Section \ref{dahmen_huang_kutyniok_lim_schwab:stabvar}, it can be shown that the resulting RBM
has stability constants as close to one as one wishes yielding for the first time  an RBM for transport and convection-diffusion problems
with this property exhibiting the same rates as the Kolmogorov widths \cite{dahmen_huang_kutyniok_lim_schwab:DePeWo13}.

In Section \ref{dahmen_huang_kutyniok_lim_schwab:RadTrans} of this report, and in \cite{dahmen_huang_kutyniok_lim_schwab:Gr13},
we present algorithms which
construct {\em explicitly} separable approximations of the form \eqref{dahmen_huang_kutyniok_lim_schwab:separation}
for the parametric transport problem of radiative transfer.
We also mention that separable approximations such as \eqref{dahmen_huang_kutyniok_lim_schwab:separation}
arise in a host of other applications;
for example, in parametric representations of PDEs with
random field input data with the aid of sparse tensor product interpolation methods;
we refer to \cite{dahmen_huang_kutyniok_lim_schwab:CoDeSc10,dahmen_huang_kutyniok_lim_schwab:ChCoDeSc13} and to the references therein.
Adaptive near-minimal rank tensor solvers for problems in high dimensional
phase space are established and analyzed in \cite{dahmen_huang_kutyniok_lim_schwab:BD}.\vspace*{-5mm}
%
\section{Anisotropic Approximations}\label{dahmen_huang_kutyniok_lim_schwab:sect2}
\vspace*{-2mm}
Let $D = (0,1)^2$ and let ${\rm curv}(\partial\Omega )$ denote the curvature of $\partial\Omega\cap D$.
The class of {\em cartoon-like functions} on $D = (0,1)^2$,
\begin{eqnarray}\label{dahmen_huang_kutyniok_lim_schwab:cartoon}
\mathcal{C}(\zeta, L,M,D) &:=& \{ f_1\chi_{\Omega } + f_2\chi_{D\setminus\Omega}: \,\,\Omega\subset D,\,   |\partial\Omega \cap D|\leq L, \partial\Omega \cap D \in C^2, \nonumber \\
&& \vspace{-70pt}  {\rm curv}(\partial\Omega )\leq \zeta,  \|f_i^{(l)}\|_{L_\infty(D)}\leq M, \, l \leq 2,\, \, i=1,2\},
\end{eqnarray}
(where the parameters $\zeta,L$ are not mutually independent)
has become a well accepted benchmark for sparse approximation in imaging \cite{dahmen_huang_kutyniok_lim_schwab:Do01}.
Compactly supported shearlet systems for $L^2(\mathbb{R}^2)$ have been introduced
in \cite{dahmen_huang_kutyniok_lim_schwab:KiKuLi12,dahmen_huang_kutyniok_lim_schwab:KuLi11} to provide (near-) optimal sparse
approximations for such classes.
We observe that such cartoons also exhibit similar features as solutions to transport problems.



{Unfortunately, even compactly supported shearlets do not comply well with quadrature  and boundary adaptation tasks
faced in variational methods for PDEs. We are therefore interested in generating locally refinable anisotropic partitions  
for which corresponding piecewise polynomial approximations realize the favorable near-optimal approximation rates for 
cartoon functions achieved by shearlet systems. Unfortunately, as shown in \cite[Chapter 9.3]{dahmen_huang_kutyniok_lim_schwab:Mi11}, 
simple triangular bisections connecting the midpoint of an edge to the opposite vertex is not sufficient for warranting such rates, 
see \cite{{dahmen_huang_kutyniok_lim_schwab:CoDyHeMi12,dahmen_huang_kutyniok_lim_schwab:ChSuXu07}} for related work. In fact, a key 
feature would be to realize a ``parabolic scaling law'' similar to the shearlet setting. By this we mean a sufficient rapid directional 
resolution by anisotropic cells whose width scales like the square of the diameter.
To achieve this we consider partitions comprised of triangles {\em and} quadrilaterals pointed out to us in \cite{dahmen_huang_kutyniok_lim_schwab:priv-com}.
We 
sketch the main ideas and refer to \cite{dahmen_huang_kutyniok_lim_schwab:DaKuLiScWe13}
for details.

Starting from some initial partition consisting of triangles and quadrilaterals, {\em refined partitions} are obtained by splitting
a given cell $Q$ of a current partition according to one of the following rules:
\begin{enumerate}
\item[(i)] Connect a vertex with the midpoint of an edge not containing the vertex.
\item[(ii)] Connect two vertices.
\item[(iii)] Connect the midpoints of two  edges which, when $Q$ is a quadrilateral, do not share any vertex.
\end{enumerate}


The types of bisections are indicated in Figure \ref{dahmen_huang_kutyniok_lim_schwab:refinement}: 

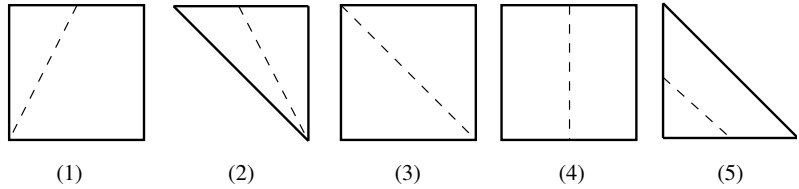
\begin{figure}[ht]
\begin{center}
\includegraphics[width=0.9\textwidth]{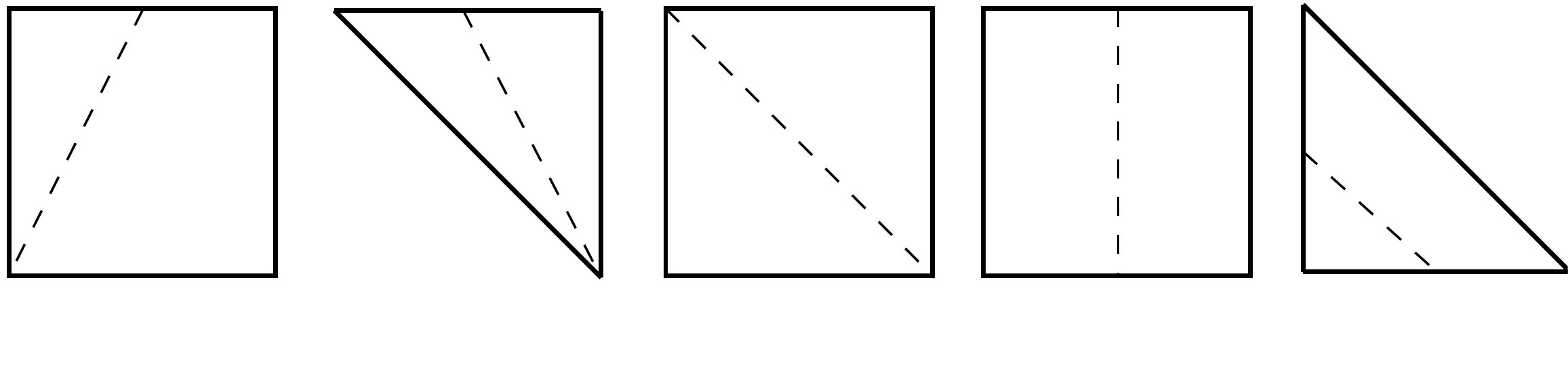}
\put(-280,5){(1)}
\put(-215,5){(2)}
\put(-152,5){(3)}
\put(-90,5){(4)}
\put(-30,5){(5)}
\end{center}
\vspace*{-3mm}
\caption{Illustration of the partion rules.}
\label{dahmen_huang_kutyniok_lim_schwab:refinement}
\end{figure}

\vspace*{-3mm}

(1), (2) are examples of (i), (3) illustrates (ii), and (4), (5) are examples for (iii). One easily checks that
these refinement rules produce only  {triangles} and quadrilaterals.
Moreover, a quadrilateral can be bisected in
$8$ possible ways whereas  a triangle can be split in $6$ possible ways. Assigning to each split type a number
in $I_Q = \{1,\dots,8\}$ when $Q$ is a quadrilateral and a number in $I_Q = \{9,\dots,14\}$ when $Q$ is a triangle,
we denote by
\begin{equation}\label{dahmen_huang_kutyniok_lim_schwab:refinement_op}
R_{\iota_Q}(Q) = \{Q_1,Q_2\} \quad \mbox{for some \,\,} \iota_Q \in I_Q,
\end{equation}
the refinement operator which replaces the cell $Q$ by its two children $Q_1,Q_2$ generated, according to the
choice $\iota_Q$, by the above split rules (i)--(iii).

For any partition $\mathcal{G}$ of $D$, let $\mathbb{P}_1( \mathcal{G}) = \{v \in L_2(D) : v|Q \in \mathbb{P}_1, Q \in \mathcal{G}\}$
be the space of piecewise affine functions on $\mathcal{G}$ and
denote by $\mathfrak{G}$ the set of all finite partitions that can be created  by successive applications of
$R_{\iota_Q}$ to define 
$$
\Sigma_N:= \bigcup \, \{\mathbb{P}_1(\mathcal{G}): \mathcal{G} \in \mathfrak{G},\, \#(\mathcal{G})\leq N\}.
$$
The next result from \cite{dahmen_huang_kutyniok_lim_schwab:DaKuLiScWe13}
shows that approximations by
elements of $\Sigma_N$ realize (and even slightly improve on)
the known rates obtained for shearlet systems for
the class of cartoon-like functions
\cite{dahmen_huang_kutyniok_lim_schwab:KuLi11}).

\begin{theorem}[\cite{dahmen_huang_kutyniok_lim_schwab:DaKuLiScWe13}]
\label{dahmen_huang_kutyniok_lim_schwab:ShearCart}
Let $f \in \mathcal{C}(\zeta,L,M,D)$ with $D = (0,1)^2$ and
assume that the discontinuity curve
$\Gamma = \partial \Omega \cap D$ is the graph of a $C^2$-function.
Then there exists a positive constant $a$ such that
\[
\inf_{\varphi \in \Sigma_N}\|f-\varphi\|_{L_2(D)} \leq {C(\zeta,L)M} \, (\log N)N^{-1},
\]
where $C(\zeta,L)$ is an absolute constant depending only on $\zeta,L$.
\end{theorem}
The proof of Theorem is based on constructing a specific sequence $\mathcal{C}_j$ of admissible partitions from $\mathfrak{G}$ where the
refinement decisions represented by $R_{\iota_Q}$ use full knowledge of the approximated function $f$. A similar sequence
of partitions is employed in Section \ref{dahmen_huang_kutyniok_lim_schwab:shear_adapt} where $\iota_Q \in I_Q$,
however, results from an {\em a posteriori} criterion described below.
We close this section by a few remarks on the structure of the $\mathcal{C}_j$.
Given $\mathcal{C}_{j-1}$, we first generate
\begin{equation}\label{dahmen_huang_kutyniok_lim_schwab:refinement_collection}
\tilde{\mathcal{C}}_{j} = \{ Q^{'} \in \tilde{R}(Q) : Q \in \mathcal{C}_{j-1} \},
\end{equation}
where $\tilde{R}$ is either $R_{\iota_Q}$ or the identity. To avoid unnecessary refinements
we define then
 $\mathcal{C}_{j}$   by replacing any pair of triangles $Q,Q^{'} \in \tilde{\mathcal{C}}_{j}$, whose union forms a parallelogram
$P$ by $P$ itself. This reduces the number of triangles in favor of parallelograms.
}
\vspace*{-3mm}

\section{Well-Conditioned Stable Variational Formulations}
\label{dahmen_huang_kutyniok_lim_schwab:stabvar}
In this section we highlight some new conceptual developments from \cite{dahmen_huang_kutyniok_lim_schwab:CoDaWe12,dahmen_huang_kutyniok_lim_schwab:DaHuScWe12,dahmen_huang_kutyniok_lim_schwab:DePeWo13} which, are, in particular, relevant for the high dimensional parametric problems addressed later below.\vspace*{-5mm}
\subsection{The General Principles}
\vspace*{-2mm}
Anisotropic structures are already exhibited by solutions of elliptic boundary value problems on polyhedral
domains  in 3D. However, related singularities are known a priori and can be dealt with by anisotropic preset
mesh refinements. Anisotropic structures of solutions to transport dominated problems can be less predictable so that
a quest for {\em adaptive anisotropic} discretization principles gains more weight. Recall that every known rigorously
founded adaptation strategy hinges in one way or the other on being able to relate a current error of an approximate solution
to the corresponding {\em residual} in a suitable norm. While classical variational formulations of elliptic problems grant exactly
such an error-residual relation, this is unclear for transport dominated problems. The first fundamental issue is therefore to {\em find}
also for such problems suitable variational formulations yielding a well conditioned error-residual relation.\vspace*{-5mm}
\subsubsection{Abstract Petrov-Galerkin Formulation}
\label{dahmen_huang_kutyniok_lim_schwab:AbsPG}
\vspace*{-3mm}
 Suppose that for a pair of Hilbert spaces $X,Y$ (with scalar products
$(\cdot,\cdot)_X, (\cdot,\cdot)_Y$ and norms $\|\cdot\|_X, \|\cdot\|_Y$),
and a given bilinear form $b(\cdot,\cdot): X\times Y$, the problem
\begin{equation}
\label{dahmen_huang_kutyniok_lim_schwab:WC1}
b(u,v)= f(v),\quad v\in Y,
\end{equation}
has for any $f\in Y'$ (the normed dual of $Y$) a unique solution $u\in X$. It is well-known that this  is equivalent to the existence of constants
$0<c_b \leq C_b<\infty$ such that
\begin{equation}
\label{dahmen_huang_kutyniok_lim_schwab:WCinfsup}
\sup_{w\in X}\sup_{v\in Y}\frac{b(w,v)}{\|w\|_X\|v\|_Y}\leq C_b,\quad \inf_{w\in X}\sup_{v\in Y}\frac{b(v,w)}{\|w\|_X\|v\|_Y}\geq c_b,
\end{equation}
and the existence of a $w\in X$ such that $b(w,v)\neq 0$ for all $v\in Y$. This means that the operator
$B:X\to Y'$, defined by $(Bu)(v):= b(u,v)$, $u\in X, v\in Y$,
is an isomorphism with condition number
$\kappa_{X,Y}(B):= \|B\|_{\mathcal{L}(X,Y')}\|B^{-1} \|_{\mathcal{L}(Y',X)}\leq C_b/c_b$.
For instance, when \eqref{dahmen_huang_kutyniok_lim_schwab:WC1}
represents a {\em convection dominated convection-diffusion problem} with
the classical choice $X=Y=H^1_0(\Omega)$, the quotient $C_b/c_b$ becomes very large.
Since
\begin{equation}
\label{dahmen_huang_kutyniok_lim_schwab:WCerres}
\|B\|_{\mathcal{L}(X,Y')}^{-1}\| Bv-f\|_{Y'}
\leq
\|u-v\|_X \leq \|B^{-1}\|_{\mathcal{L}(Y',X)}\| Bv-f\|_{Y'},
 \end{equation}
the error $\|u-v\|_X$ can then {\em not} be tightly estimated by the residual $\| Bv-f\|_{Y'}$.\vspace*{-3mm}

\subsubsection{Renormation}
\label{dahmen_huang_kutyniok_lim_schwab:Renorm}
\vspace*{-3mm}
On an abstract level the following principle has
surfaced in a number of different contexts such as
{\em least squares methods} (see e.g. \cite{dahmen_huang_kutyniok_lim_schwab:BoGu09})
and the so-called, more recently emerged
{\em Discontinuous Petrov Galerkin (DPG) methods},
see e.g.
\cite{
dahmen_huang_kutyniok_lim_schwab:BaMo84,dahmen_huang_kutyniok_lim_schwab:DaHuScWe12,
dahmen_huang_kutyniok_lim_schwab:DeGo10,dahmen_huang_kutyniok_lim_schwab:DeGo11}
and the references therein.
The idea is to fix a norm, $\|\cdot\|_Y$, say, and modify the norm for
$X$ so that the corresponding operator even becomes an {\em isometry}.
More precisely, define\vspace*{-1mm}
\begin{equation}
\label{dahmen_huang_kutyniok_lim_schwab:CWhatX}
\|u\|_{\hat X}:= \sup_{v\in Y}\frac{b(u,v)}{\|v\|_Y}= \|B u\|_{Y'} = \|R_Y^{-1}Bu\|_Y,
\end{equation}
where $R_Y: Y\to Y'$ is the Riesz map defined by $(v,z)_Y= (R_Y v)(z)$.
The following fact is readily verified, see e.g. \cite{dahmen_huang_kutyniok_lim_schwab:DaHuScWe12,dahmen_huang_kutyniok_lim_schwab:We13}.

\begin{remark}
\label{dahmen_huang_kutyniok_lim_schwab:WCrem1}
One has $\kappa_{\hat X,Y}(B)=1$, i.e., \eqref{dahmen_huang_kutyniok_lim_schwab:WCinfsup} holds with $c_b=C_b=1$ when $\|\cdot\|_X$ is replaced by $\|\cdot\|_{\hat X}$.

 Alternatively, fixing $X$ and redefining $\|\cdot\|_Y$ by $\|v\|_{\hat Y}:= \|B^*v\|_{X'}$, one has $\kappa_{X,\hat Y}(B)=1$,  see \cite{dahmen_huang_kutyniok_lim_schwab:DaHuScWe12}.
Both possibilities lead to the error residual relations
\begin{equation}
\label{dahmen_huang_kutyniok_lim_schwab:WCperfect}
\|u - w\|_X = \|f - Bw\|_{\hat Y'},\quad \|u - w\|_{\hat X} = \|f - Bw\|_{Y'}, \quad u, w\in X.
\end{equation}
\end{remark}
\vspace*{-5mm}
\subsection{Transport Equations}\label{dahmen_huang_kutyniok_lim_schwab:transport}
\vspace*{-3mm}
Several variants of these principles are applied and analyzed in detail in
\cite{dahmen_huang_kutyniok_lim_schwab:CoDaWe12} for convection-diffusion equations.
We concentrate
in what follows on the limit case for vanishing viscosity, namely pure transport equations.
For simplicity we consider the domain $D= (0,1)^d$, $d= 1,2,3$, with
$\Gamma := \partial D$, denoting as usual by $\vec{n}=\vec{n}(x)$
the unit outward normal at $x\in \Gamma$ (excluding the four corners, of course).
Moreover, we consider
velocity fields $\vec{b}(x)$, $x\in D$, which for simplicity will always
be assumed to be differentiable, i.e., $\vec{b}(x) \in C^1(\overline{D})^d$.
Likewise $c(x) \in C^0(\overline{D})$ will serve as the reaction
term in the {\em first order transport equation} \vspace*{-1mm}
\begin{equation}
 \vec{b} \cdot \nabla u + cu  = \;\mbox{$f_\circ$ in $D$}\,, \label{dahmen_huang_kutyniok_lim_schwab:2.1}
 \quad
u  = \; \mbox{$g$ on $\Gamma_-$}\,,
\end{equation}
where
$\Gamma_\pm :=  \{x \in \partial D: \; \pm \vec{b} (x) \cdot \vec{n}(x) > 0\}$
denotes the {\em inflow, outflow boundary}, respectively.
%
Furthermore, to simplify the exposition we shall always assume
that 
%
$2c -  
 \nabla\cdot  \vec{b} \geq c_0 > 0$ in $D$
%
holds.

A priori there does not seem to be any  ``natural'' variational formulation. Nevertheless, the above principle
can be invoked as follows.
Following e.g. \cite{dahmen_huang_kutyniok_lim_schwab:DaHuScWe12},
one can show that  the associated bilinear form with derivatives on the test functions
\begin{equation}
\label{dahmen_huang_kutyniok_lim_schwab:auv}
b(w,v):=   \displaystyle\int_D \,w(- \vec{b} \cdot {\nabla} v + v (c- {\nabla} \cdot \vec{b}))  \;dx,
\end{equation}
%
is trivially bounded on $L_2(D)\times W_0(-\vec{b},D)$, where
\begin{equation}
\label{dahmen_huang_kutyniok_lim_schwab:W0}
W_0(\mp \vec{b},D)
:=
{\rm clos}_{\|\cdot\|_{W(\vec{b},D)}}\{ v\in C^1(D)\cap C(\overline{D}),\,
v\mid_{\Gamma_{\pm}}\equiv 0\}
\end{equation}
and
\begin{equation}
\label{dahmen_huang_kutyniok_lim_schwab:graph}
\|v\|_{W(\vec{b},D)}:=
\left( \|v\|_{L_2(D)}^2 + \int_D |\vec{b}\cdot\nabla v|^2\,dx\right)^{1/2}.
\end{equation}
Moreover,  the trace $\gamma_-(v)$ on the inflow boundary
exists and is contained in
$L_2(\Gamma_-, {|\vec{b}\cdot\vec{n}|})$ for $v\in W_0(\vec{b},D)$, endowed with the norm
$\|g\|^2_{L_2(\Gamma_\pm, {|\vec{b}\cdot\vec{n}|})} = \int_{\Gamma_\pm} |g|^2| {\vec{b}\cdot\vec{n}|} ds$
so that
\begin{equation}
\label{dahmen_huang_kutyniok_lim_schwab:ell}
f(v)
:= (f_\circ,v)
+ \int_{\Gamma_-}g\gamma_-(v)|\vec{b}\cdot\vec{n}|ds
\end{equation}
belongs to  $ (W_0(\vec{b},D))'$ and the variational problem
\begin{equation}
\label{dahmen_huang_kutyniok_lim_schwab:varprob}
b(u,v) = f(v),\quad v\in W_0(-\vec{b},D)
\end{equation}
possesses a unique solution in $L_2(D)$ which, when regular enough,
coincides with the classical solution of \eqref{dahmen_huang_kutyniok_lim_schwab:2.1},
see  \cite[Theorem 2.2]{dahmen_huang_kutyniok_lim_schwab:DaHuScWe12}.

Moreover, since $X=L_2(D)=X'$,
the quantity $\|v\|_Y:= \|B^* v\|_{L_2(D)}$ is an equivalent norm on $W_0(-\vec{b},D)$,
see \cite{dahmen_huang_kutyniok_lim_schwab:DaHuScWe12},
and Remark \ref{dahmen_huang_kutyniok_lim_schwab:WCrem1}
applies, 
%
 i.e.,
\begin{equation}
\label{dahmen_huang_kutyniok_lim_schwab:isos}
\|B\|_{\mathcal{L}(L_2(D),(W_0(\vec{b},D))')} = \|B^*\|_{\mathcal{L}(W_0(\vec{b},D),L_2(D))} =1,
\end{equation}
see \cite[Proposition 4.1]{dahmen_huang_kutyniok_lim_schwab:DaHuScWe12}.
One could also reverse the roles of test and trial space (with the inflow boundary conditions
being then essential ones) but the present formulation imposes least regularity on the solution which will be essential in the next section.
Note that whenever a PDE is written as a first order system, $X$ can always be arranged as an $L_2$-space.

Our particular interest concerns the {\em parametric case}, i.e., the constant convection field $\vec{s}$ in
\begin{equation}
\label{dahmen_huang_kutyniok_lim_schwab:rad1.1}
\begin{array}{rcl}
  \vec{s}\cdot \nabla u(x,\vec{s})+\kappa (x)u(x,\vec{s})&=&
f_\circ(x),\quad x\in D\subset \mathbb{R}^d,\,\, d=2,3,
\\
u(x,\vec{s})&=& g(x,\vec{s}),\, x\in \Gamma_-(\vec{s}),
\end{array}
\end{equation}
may {\em vary} over  a set of directions $\mathcal{S}$ so that
  now the solution $u$ also depends on the transport direction $\vec{s}$.
In \eqref{dahmen_huang_kutyniok_lim_schwab:rad1.1}
 {and the following} we assume that
${\rm ess}\inf_{x\in D} \kappa(x) \geq 0$.
Thus, for instance, when $\mathcal{S}=S^2$, the unit $2-$sphere,
$u$ is considered as a function of five variables,
namely $d=3$ spatial variables and parameters from a two-dimensional set $\mathcal{S}$.
This is the simplest example of a kinetic equation forming a core constituent in {\em radiative transfer} models.
The  in- and outflow boundaries now depend on $\vec{s}$:
\begin{equation}\label{dahmen_huang_kutyniok_lim_schwab:Gammapm}
\Gamma_\pm(\vec{s}):= \{x\in \partial D: \mp \vec{s}\cdot {\bf n}(x)<0\},\qquad \vec{s}\in \mathcal{S}\;.
\end{equation}
Along similar lines one can determine $u$ as a function of $x$ and $\vec{s}$ in $X= L_2(D\times \mathcal{S})$
as the solution of a variational problem with test space $Y 
 :={\rm clos}_{\|\cdot\|_{W(D\times \mathcal{S})}} \{v\in C(\mathcal{S},C^1(D)): v|_{\Gamma_\pm}\equiv 0\}$ with
 $  \|v\|_{W(D\times \mathcal{S})}^2 :=
\|v\|_{L_2(D\times \mathcal{S})}^{2} + \int_{\mathcal{S}\times D}|\vec{s}\cdot \nabla v|^2 dx d\vec{s}$. Again this formulation
requires minimum regularity. Since later we shall discuss yet another formulation,
imposing stronger regularity conditions, we refer to \cite{dahmen_huang_kutyniok_lim_schwab:DaHuScWe12} for details.\vspace*{-7mm}
\subsection{$\delta$-Proximality and Mixed Formulations}
\vspace*{-3mm}
It is initially not clear how to exploit
\eqref{dahmen_huang_kutyniok_lim_schwab:WCperfect} numerically
since the perfect inf-sup stability on the infinite
dimensional level is {\em not} automatically inherited by finite dimensional subspaces
$X_h\subset X, Y_h\subset Y$ of equal dimension.
However, given $X_h\subset \hat X$, one can identify the
``ideal'' test space $Y(X_h)=R_Y^{-1}B(X_h)$
which may be termed ideal because
\begin{equation}
\label{dahmen_huang_kutyniok_lim_schwab:WCinfsuph}
\sup_{w\in X_h}\sup_{v\in Y(X_h)}\frac{b(w,v)}{\|w\|_X\|v\|_Y}=   \inf_{w\in X_h}\sup_{v\in Y(X_h)}\frac{b(v,w)}{\|w\|_X\|v\|_Y}=1,
\end{equation}
 see \cite{dahmen_huang_kutyniok_lim_schwab:DaHuScWe12}. In particular, this means that the solution
$u_h\in X_h$ of the corresponding {\em Petrov-Galerkin} scheme
\begin{equation}
\label{dahmen_huang_kutyniok_lim_schwab:WCidPG}
b(u_h,v)= f(v),\quad v\in Y(X_h),
\end{equation}
realizes the best $\hat X$-approximation to the solution $u$ of \eqref{dahmen_huang_kutyniok_lim_schwab:WC1}, i.e.,
\begin{equation}
\label{dahmen_huang_kutyniok_lim_schwab:WCbest}
\|u-u_h\|_{\hat X}= \inf_{w\in X_h}\|u-w\|_{\hat X}.
\end{equation}
Of course, unless $Y$ is an $L_2$ space, the ideal test space $Y(X_h)$  is,
in general, not computable exactly.
To retain stability it is natural to look for a numerically computable test space $Y_h$ that
is sufficiently close to $Y(X_h)$.

One can pursue several different strategies to obtain numerically
feasible test spaces $Y_h$.
When   \eqref{dahmen_huang_kutyniok_lim_schwab:WC1} is a discontinous Galerkin formulation
one can choose $Y$ as a product space over the given
partition, again with norms induced by the graph norm for the
adjoint $B^*$ so that the approximate inversion of the Riesz map $R_Y$
can be {\em localized} \cite{dahmen_huang_kutyniok_lim_schwab:DeGo10,dahmen_huang_kutyniok_lim_schwab:DeGo11}.
An alternative, suggested in \cite{dahmen_huang_kutyniok_lim_schwab:CoDaWe12,dahmen_huang_kutyniok_lim_schwab:DaHuScWe12}, is based on noting that by \eqref{dahmen_huang_kutyniok_lim_schwab:CWhatX}
the ideal Petrov Galerkin solution $u_h$ from \eqref{dahmen_huang_kutyniok_lim_schwab:WCidPG} is a
{\em minimum residual} solution in $Y'$, i.e., $u_h= {\rm argmin}_{w\in X_h}
\|f- Bw\|_{Y'}$ whose normal equations read $(f- B u_h,Bw)_{Y'}=0$,
$w\in X_h$. Since the inner product $(\cdot,\cdot)_{Y'}$
is numerically hard to access, one can  write  $(f- B u_h,Bw)_{Y'}= \langle R_Y^{-1}(f- Bu_h),Bw\rangle$,
where the dual pairing $\langle \cdot,\cdot\rangle$ is now induced by  the standard $L_2$-inner product.
Introducing as an auxiliary variable the ``lifted residual''
\begin{equation}
\label{dahmen_huang_kutyniok_lim_schwab:WClift}
y= R_Y^{-1}(f- Bu_h),
\end{equation}
or equivalently $(R_Y y)( v) =\langle R_Y y,v\rangle = (y,v)_Y = \langle f-  Bu_h,v\rangle$, $v\in Y$, one can show that \eqref{dahmen_huang_kutyniok_lim_schwab:WCidPG} is equivalent to the saddle point problem
\begin{equation}
\label{dahmen_huang_kutyniok_lim_schwab:system-full}
\begin{array}{lccl}
  \langle  R_{Y} y,v\rangle + b (u_h,v)
   & = & \langle f,v\rangle, & v\in Y,\\
b(w, y)
&=& 0,& w \in X_h,
\end{array}
\end{equation}
which involves only standard $L_2$-inner products, see \cite{dahmen_huang_kutyniok_lim_schwab:DaHuScWe12,dahmen_huang_kutyniok_lim_schwab:DePeWo13}.
\begin{remark}
\label{dahmen_huang_kutyniok_lim_schwab:rem2}
When working with $X, \hat Y$ instead of $\hat X, Y$, one has $R_Y = B R_X^{-1}B^*$ and hence, when
 $X= L_2(D)$ as in \eqref{dahmen_huang_kutyniok_lim_schwab:auv}, one has $R_Y = BB^*$ (see also \cite{dahmen_huang_kutyniok_lim_schwab:CaMaMcRu01}).
 \end{remark}
Since the test space $Y$ is still infinite dimensional, a numerical realization would require finding a (possibly small) subspace $V\subset Y$
such that the analogous saddle point problem with $Y$ replaced by $V$ is still inf-sup stable. The relevant condition on $V$
can be described by the notion of $\delta$-{\em proximality}
introduced in \cite{dahmen_huang_kutyniok_lim_schwab:DaHuScWe12},
see also \cite{dahmen_huang_kutyniok_lim_schwab:CoDaWe12}.
We  recall the formulation from \cite{dahmen_huang_kutyniok_lim_schwab:DePeWo13}:
\emph{$V\subset Y$ is $\delta$-proximal for $X_h\subset \hat X$ if},
for some $\delta \in (0,1)$,
with
$P_{Y,V}$ denoting the $Y$-orthogonal projection from $Y$ to $V$,
\begin{equation}
\label{dahmen_huang_kutyniok_lim_schwab:delta-prox}
\|(I- P_{Y,V})R_{Y}^{-1}B  w\|_{Y} \leq \delta \|R_{Y}^{-1}Bw\|_{Y},
\quad w\in X_h
\;.
\end{equation}
%
%
\begin{theorem}
{\rm \cite{dahmen_huang_kutyniok_lim_schwab:CoDaWe12,dahmen_huang_kutyniok_lim_schwab:DaHuScWe12,dahmen_huang_kutyniok_lim_schwab:DePeWo13}}
\label{dahmen_huang_kutyniok_lim_schwab:stab2}
Assume that for given $X_h\times V \subset X\times Y$ the test space
$V$ is $\delta$-proximal for $X_h$, i.e. \eqref{dahmen_huang_kutyniok_lim_schwab:delta-prox} is satisfied.
Then,
the solution $ (u_{X_h,V}, y_{X_h,V})\in X_h \times V$
of the saddle point problem
\begin{equation}
\label{dahmen_huang_kutyniok_lim_schwab:system-disc}
\begin{array}{lccl}
  \langle R_{Y} y_{X_h,V},v\rangle + b (u_{X_h,V},v)
   & = & \langle f,v\rangle, & v\in V,\\[1mm]
b(w,y_{X_h,V})
&=& 0,& w \in X_h,
\end{array}
\end{equation}
satisfies
\begin{equation}
\label{dahmen_huang_kutyniok_lim_schwab:BAP1}
\|u- u_{X_h,V} \|_{\hat X}\leq \frac{1}{1-\delta}\inf_{w\in X_h}\|u-w\|_{\hat X}.
\end{equation}
and
\begin{equation}
\label{dahmen_huang_kutyniok_lim_schwab:BAP2}
\|u- u_{X_h,V} \|_{\hat X}+ \|y- y_{X_h,V}\|_{Y}\leq \frac{2 }{1-\delta}
\inf_{w\in X_h}\|u-w\|_{\hat X}.
\end{equation}
Moreover, one has
\begin{equation}
\label{dahmen_huang_kutyniok_lim_schwab:inf-sup-0}
\inf_{w\in X_h} \sup_{v\in V}\frac{b (w, v)}{\|v\|_{Y}\|q\|_{\hat X}} \geq  {\sqrt{1-\delta^2}} .
\end{equation}
Finally, \eqref{dahmen_huang_kutyniok_lim_schwab:system-disc} is equivalent to the Petrov-Galerkin scheme
\begin{equation}
\label{dahmen_huang_kutyniok_lim_schwab:WCnearbestPG}
b(u_{X_h,V},v) = f(v),\quad v\in Y_h := P_{Y,V}(R_Y^{-1}B(X_h))= P_{Y,V}(Y(X_h)).
\end{equation}
\end{theorem}
The central message is that the Petrov-Galerkin scheme \eqref{dahmen_huang_kutyniok_lim_schwab:WCnearbestPG} can be realized {\em without} computing
a basis for the test space $Y_h$, which for {\em each} basis function could require solving a problem of the size ${\rm dim}\,V$,
by solving instead the saddle point problem \eqref{dahmen_huang_kutyniok_lim_schwab:system-disc}. Moreover, the stability of both problems is goverend by
the $\delta$-proximality of $V$. As a by-product, in view of \eqref{dahmen_huang_kutyniok_lim_schwab:WClift}, the solution component $y_{X_h,V}$ approximates
the exact lifted residual $R_Y^{-1}(f- B u_{X_h,V})$ and, as pointed out below, can be used for an a posteriori error control.

The problem \eqref{dahmen_huang_kutyniok_lim_schwab:system-disc}, in turn, can be solved with the aid of an {\em Uzawa iteration} whose efficiency relies again
on $\delta$-proximality.
For $k = 0,\dots$,  
solve
\begin{eqnarray}
\label{dahmen_huang_kutyniok_lim_schwab:uzawa}
\langle R_Y  y^k,  v\rangle &=& \langle f-B u^k,v\rangle, \quad v \in V, \nonumber\\
(u^{k+1},w)_{\hat X} &=& (u^k,w)_{\hat X}+\langle  B^*  y^k, w\rangle, \quad w \in X_h.
\end{eqnarray}
Thus, each iteration requires solving a symmetric positive definite Galerkin problem in $V$
for the approximate lifted residual.

\begin{theorem}
\label{dahmen_huang_kutyniok_lim_schwab:th7.1}
  Assume that \eqref{dahmen_huang_kutyniok_lim_schwab:delta-prox} is satisfied. 
Then the iterates
generated by the  scheme \eqref{dahmen_huang_kutyniok_lim_schwab:uzawa} converge to $u_{X_h,V}$ and
\begin{equation}
\label{dahmen_huang_kutyniok_lim_schwab:uhdeltauk}
\|u_{X_h,V}-u^{k+1}\|_{\hat X} \leq \delta \|u_{X_h,V}-u^{k}\|_{\hat X},\quad k=0,1,2,\ldots.
\end{equation}
\end{theorem}
\subsection{Adaptive Petrov-Galerkin Solvers on Anisotropic Approximation Spaces}
\vspace*{-3mm}
The benefit of the above saddle point formulation is not only that it saves us the explicit calculation of the test basis functions
but that it provides also an {\em error estimator} based on the {\em lifted residual} $y_h=y_h(u_{X_h,V},f) $ defined by  the first row of \eqref{dahmen_huang_kutyniok_lim_schwab:system-disc}.
\vspace*{-5mm}
\subsubsection{Abstract $\delta$-Proximinal Iteration}
\label{dahmen_huang_kutyniok_lim_schwab:AbsIt}
\vspace*{-3mm}
In fact, it is shown in \cite{dahmen_huang_kutyniok_lim_schwab:DaHuScWe12} that when
$V_h\subset Y$ is even $\delta$-proximal for $X_h + B^{-1}F_h$,
with some finite dimensional  subspace $F_h\subset Y'$,
one has
\begin{equation}
\label{dahmen_huang_kutyniok_lim_schwab:rh}
(1-\delta)\|f_h - B w\|_{Y'} \leq \|y_h(w,f_h)\|_Y \leq \|f_h - B w\|_{Y'} ,\quad w\in X_h,
\end{equation}
where $f_h \in F_h$ is an approximation of $f \in Y'$.
The space $F_h$ controls which components of $f$ are accounted for in the error estimator.
The  {term} $f-f_h$ is a {\em data  oscillation} error  as
encountered in adaptive finite element methods.
It follows that the current error of the Petrov-Galerkin approximation $u_{X_h,V}$
is  controlled from below and above by the quantity $\|y_h\|_Y$.
This can be used to formulate the {\em adaptive} Algorithm \ref{dahmen_huang_kutyniok_lim_schwab:adapt}
that can be proven to give rise to a {\em fixed error reduction} per step.
Its precise formulation can be found in \cite[\S~4.2]{dahmen_huang_kutyniok_lim_schwab:DaHuScWe12}.
It is shown in \cite[Proposition 4.7]{dahmen_huang_kutyniok_lim_schwab:DaHuScWe12} that Algorithm \ref{dahmen_huang_kutyniok_lim_schwab:adapt} below terminates after finitely many steps and outputs an approximate
solution $\bar u$ satisfying $\|u-\bar u\|_{\hat X} \leq \epsilon$.\vspace*{-4mm}

\begin{algorithm}[htb]
  \caption{adaptive algorithm}
  \label{dahmen_huang_kutyniok_lim_schwab:adapt}
  \begin{algorithmic}[1]
%
\State Set  
target accuracy $\epsilon$, initial guess $\bar u=0$, initial error bound $e=\|f\|_{Y'}$,
parameters $\rho,\eta, \alpha_1,\alpha_2\in (0,1)$,
initial trial and $\delta$-proximal test spaces $X_h, V_h$;
\While{$e> \epsilon$} solve \eqref{dahmen_huang_kutyniok_lim_schwab:system-disc}
within accuracy $\alpha_1\rho$ (e.g. by an Uzawa iteration with initial guess $\bar u$)
to obtain an approximate solution pair  $(\hat y, \hat u) \in V_h\times X_h$;
\State 
enlarge $X_h$ to $X_{h,+}$ in such a way that
\begin{equation}
\label{dahmen_huang_kutyniok_lim_schwab:expand}
\inf_{g\in X_{h,+}}\|B^*\hat y - g\|_{\hat X'} \leq \eta \|B^* {\hat y}\|_{\hat X'}\quad\mbox{and set}\quad r := {\rm argmin}_{g\in X_{h,+}} \|B^*{\hat y} - g\|_{\hat X'};
\end{equation}
\State
compute $X_{h'}\supset X_h, F_{h'}\supset F_h$, $f_h\in BX_{h'} + F_{h'}$ such that
$\|f- f_h\|_{Y'}\leq \alpha_2 \rho e$;
\State
set $X_h + X_{h,+} + X_{h'} \to X_h$, $\rho e \to e$, and choose a $\delta$-proximal subspace $V_h$ for $X_h$;\\
set $\hat u + r_X \to \bar u$.
 \EndWhile
  \end{algorithmic}
\end{algorithm}
\vspace*{-3mm}

\subsubsection{Application to Transport Equations}\label{dahmen_huang_kutyniok_lim_schwab:shear_adapt}
\vspace*{-3mm}
We adhere to the setting described in Section \ref{dahmen_huang_kutyniok_lim_schwab:transport}, i.e., $X=\hat X = L_2(D)$,  $\hat Y = Y = W_0(-\vec{b},D)$,
and $R_Y= BB^*$.

The trial spaces  that we now denote by $X_j$ to emphasize the nested construction below, are spanned by discontinuous piecewise linear functions on a mesh composed of cells from collections  {$\mathcal{C}_j$}, i.e.,
\begin{equation}
\label{dahmen_huang_kutyniok_lim_schwab:Xj}
X_{ {j}} = \mathbb{P}_1(\mathcal{C}_{{j}}),\quad j\geq 0,
\end{equation}
where the collections $\mathcal{C}_{{j}}$ are derived from collections $\tilde{\mathcal{C}}_{{j}}$ of the type
\eqref{dahmen_huang_kutyniok_lim_schwab:refinement_collection} as described in Section \ref{dahmen_huang_kutyniok_lim_schwab:sect2}.

Given $X_{{j}}$ of the form \eqref{dahmen_huang_kutyniok_lim_schwab:Xj}, the test spaces $V_{{j}}$ are defined by
\begin{equation}
\label{dahmen_huang_kutyniok_lim_schwab:Zj}
V_{{j}}:= \mathbb{P}_2({\mathcal{G}}_{{j}})\cap C(D) \quad \mbox{with}\quad {\mathcal{G}}_{{j}} := \{R^{iso}(Q): Q \in \mathcal{C}_{{j}}\},
\;
\end{equation}
where $R^{iso}(Q) = \{Q \cap P_i : i = 1,\dots,4\}$
is defined as follows. Let $P$ be a parallelogram containing $Q$ and sharing
at least three vertices with $Q$.
(There exist at most two such parallelograms and  we choose one of them).
Then the parallelograms $P_i$  result from a  dyadic refinement of $P$.
As pointed out later, the test spaces
$V_{{j}}$ constructed in this way, appear to be
sufficiently large to ensure $\delta$-proximality
for $X_{{j}}$ for $\delta$ significantly smaller than one
uniformly with respect to $j$.

Since the test spaces $V_j$ are determined by the trial spaces $X_j$,
the crucial step is to generate $X_{j+1}$ by enlarging $X_j$
based on an a posteriori criterion that ``senses''  directional information.
This, in turn, is tantamount to a possibly
anisotropic adaptive refinement of  $\mathcal{C}_j$
leading to the updated spaces for the next iteration sweep of the form
\eqref{dahmen_huang_kutyniok_lim_schwab:uzawa}.
The idea is to use a greedy strategy based on the largest ``fluctuation coefficients''.
To describe this, we denote  for each $\iota_Q \in I_Q$ by $\Psi_{R_{\iota_Q}(Q)}$  an orthonormal wavelet type basis for the difference space
$\mathbb{P}_1(R_{\iota_Q}(Q)) \ominus \mathbb{P}_1(Q)$.
We then set
\begin{equation}
\label{dahmen_huang_kutyniok_lim_schwab:two-level-basis}
 \Psi_j = \{\psi_{\gamma} \in\Psi_{R(Q)} : Q \in \mathcal{C}_{j-1}\},
\end{equation}
where
$
\Psi_{R(Q)} = \bigcup_{\iota_Q \in I_Q} \Psi_{R_{\iota_Q}(Q)}
$.
Initializing ${\mathcal{C}_0}$ as a uniform partition (on a low level),
we define for some fixed $\theta \in (0,1)$
\[
 T_j = \theta \cdot \max_{\psi_{\gamma} \in \Psi_j} |\langle B^*r^K_j,\psi_{\gamma}\rangle|
\]
for $j > 0$, where $\Psi_j$ is the two level basis defined in \eqref{dahmen_huang_kutyniok_lim_schwab:two-level-basis}
and $r^K_j =y^K$ is the lifted residual from the first row of the Uzawa iteration.
Then, for each $Q \in \mathcal{C}_{j-1}$, we define  its refinement $\tilde{R}(Q)$
(see the remarks following  \eqref{dahmen_huang_kutyniok_lim_schwab:refinement_collection})
by
\[
\tilde{R}(Q) := \left\{
\begin{array}{rl}
\{Q\}, & \quad \mbox{if} \quad \max_{\psi_{\gamma} \in \Psi_{R(Q)}} |\langle B^*r^K_j,\psi_{\gamma}\rangle| \leq T_j,\\
R_{\hat{\iota}_Q}(Q), & \quad \mbox{otherwise},
\end{array}
\right.
\,\,
\]
where $\hat{\iota}_{Q}$ is chosen to maximize
$\max_{\psi_{\gamma} \in \Psi_{R_{\iota_Q}(Q)}} |\langle B^*r^K_j,\psi_{\gamma}\rangle|$ among all $\iota_Q \in I_Q$.
One can then check whether this enrichment yields a sufficiently accurate
$L_2$-approximation of $B^*r^K_j$
(step 3 of Algorithm \ref{dahmen_huang_kutyniok_lim_schwab:adapt}).
{In this case, we adopt} $\mathcal{C}_j$.
Otherwise, the procedure is repeated for a smaller threshold $\theta$.

\subsection{Numerical Results}
\vspace*{-3mm}

We provide some numerical experiments to  {illustrate} the performance of the previously introduced anisotropic adaptive
scheme for first order linear transport equations and
refer to \cite{dahmen_huang_kutyniok_lim_schwab:DaKuLiScWe13} for further tests.
We monitor $\delta$-proximality by computing
\begin{equation}\label{eq:DeltaProx}
\frac{\inf_{\phi \in V_j}\|u_j -u_j^K - B^*\phi\|_{L_2([0,1]^2)}}{\|u_j-u_j^K\|_{L_2([0,1]^2)}},
\end{equation}
where $u_j = {\rm argmin}_{v_j \in X_j}\|u-v_j\|_{L_2(D)}$.
This is only a lower bound of the $\delta$-proximality constant $\delta$
for one  particular choice of $w$ in
\eqref{dahmen_huang_kutyniok_lim_schwab:delta-prox} which
 coincides with the choice of $w$ in the proof
in \cite{dahmen_huang_kutyniok_lim_schwab:DaHuScWe12}.
In the  following experiment, the number $K$ of Uzawa iterations
is for simplicity set to $K=10$.
One could as well employ an early termination of the inner iteration
based on a posteriori control of the lifted residuals $r^k_j$.

 We  consider the transport equation \eqref{dahmen_huang_kutyniok_lim_schwab:2.1} with zero boundary condition $g=0$,  convection field $\vec{b} =(x_2,1)^T$,
and right hand side $ f = \chi_{\{x_1 > x_2^2/2\}} + 1/2 \cdot \chi_{\{x_1 \leq x_2^2/2\}}$  so that the  solution exhibits a discontinuity along the
curvilinear shear layer given by $x_1 = \frac12 x_2^2$.

\begin{figure}[ht]
\begin{center}
\includegraphics[width=0.98\textwidth]{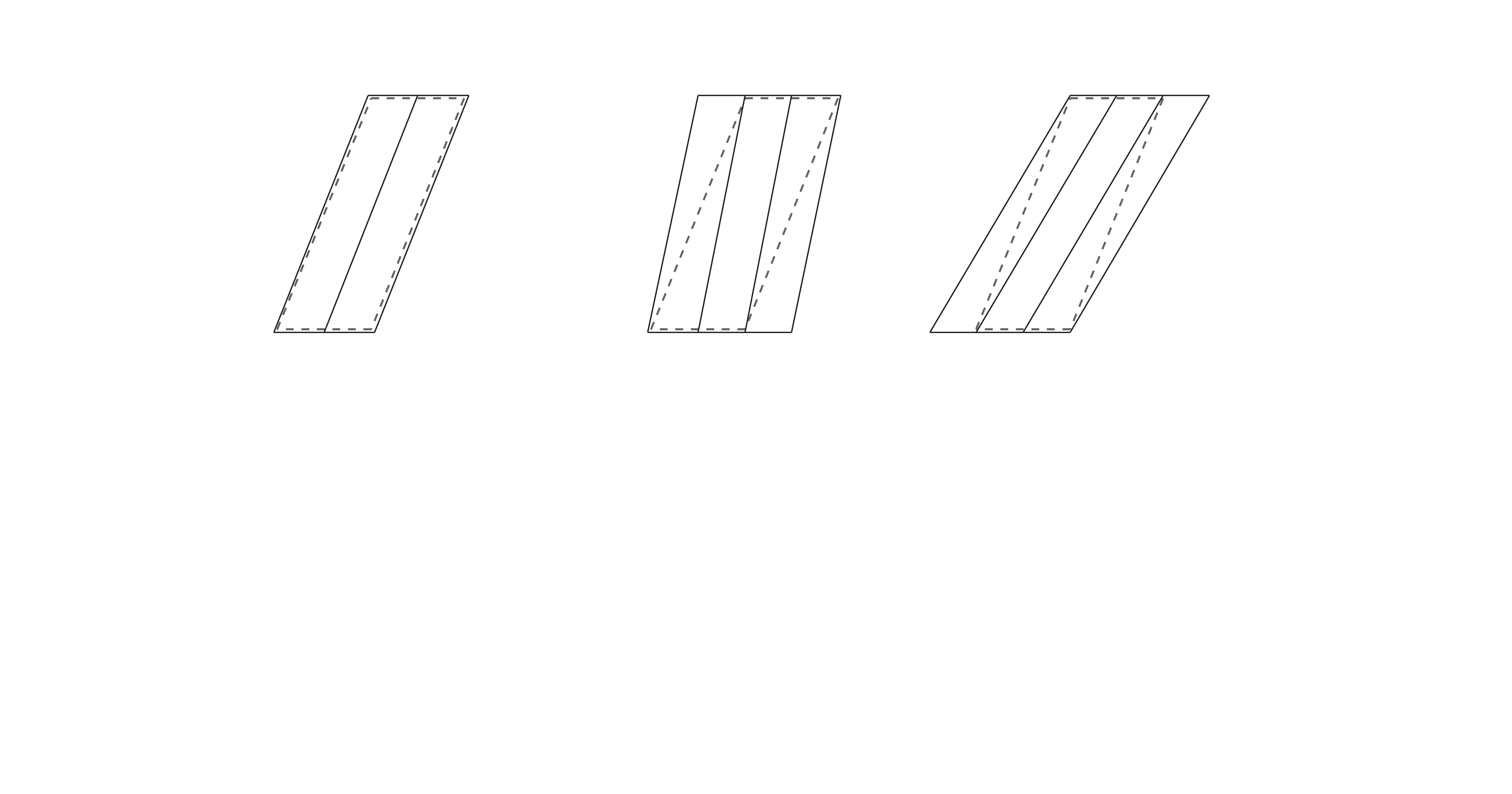} 
\put(-250,80){(a)}
\put(-170,80){(b)}
\put(-90,80){(c)}
\end{center}
\vspace*{-30mm}
\caption{Possible directional adjustments are illustrated for a  parallelogram $P$ (dashed line). (a): rule (iii) of
Section \ref{dahmen_huang_kutyniok_lim_schwab:sect2} yields two parallelograms with the same ``direction''. (b), (c): applying rule (i) twice,
changes the anisotropic direction slightly. The three refined parallelograms depicted in (b), (c) illustrate the results of a possible merging
of adjacent triangles.}
\label{dahmen_huang_kutyniok_lim_schwab:refine_implementation}
\end{figure}

In this numerical example we actually explore ways of reducing the relatively
large number of possible splits corresponding to
the operators $R_{\iota_Q}$, $\iota_Q\in I_Q$,
while still realizing the parabolic scaling law.
In fact, we confined the cells to intersections
of parallelograms $P$ and their intersections with the domain $D$,
much in the spirit of shearlet systems, employing
anisotropic refinements as illustrated in Figure
\ref{dahmen_huang_kutyniok_lim_schwab:refine_implementation}
as well as the isotropic refinement $R^{iso}$.
Permitting occasional overlaps of parallelograms,
one can even avoid any interior triangles, apparently without degrading the
accuracy of the adaptive approximation.
The  general refinement scheme described in Section
\ref{dahmen_huang_kutyniok_lim_schwab:sect2} 
covers the presently proposed one as a special case, except, of course, for the possible overlap of cells.
\vspace*{-10mm}

\begin{figure}[h]
\begin{center}
\includegraphics[height=3in]{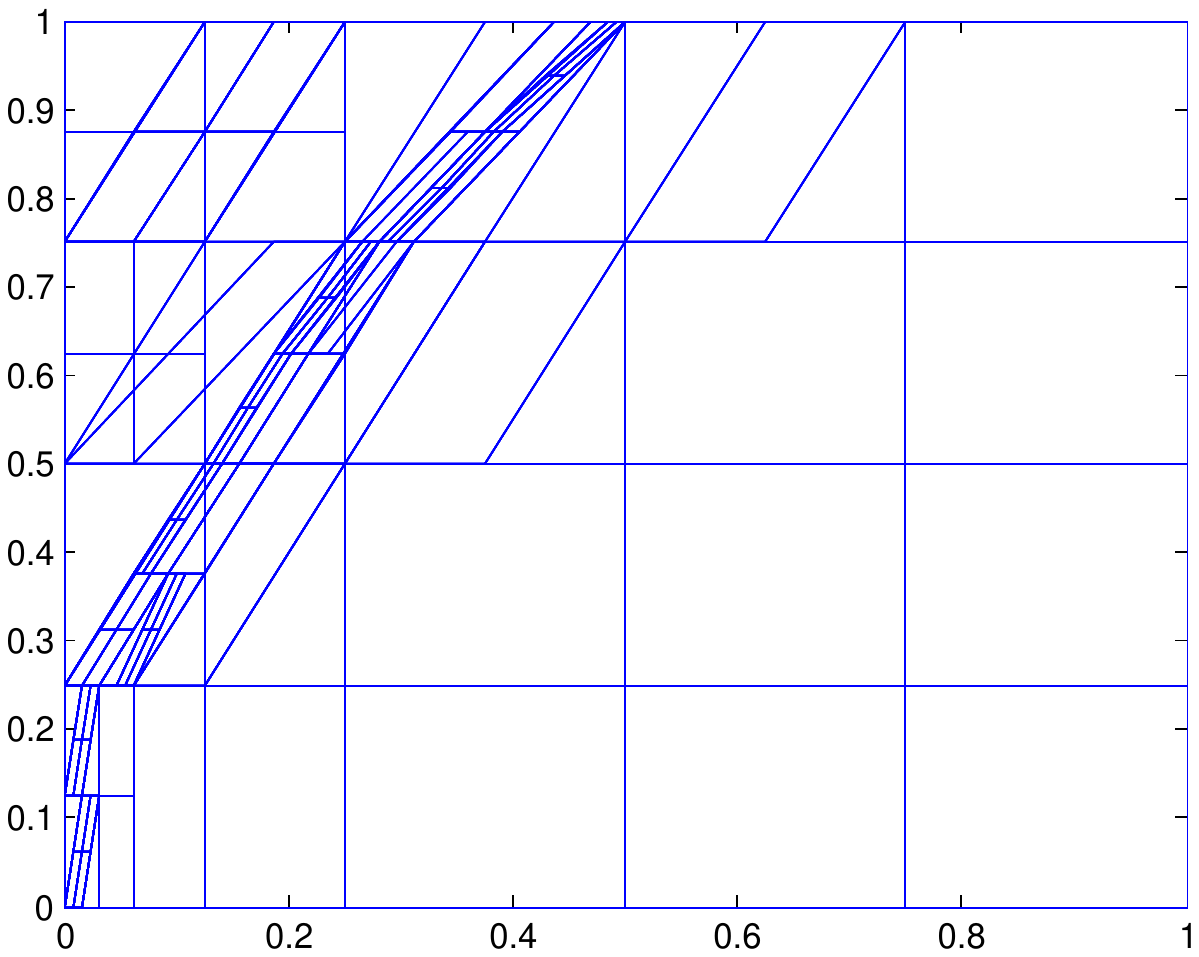}\hspace{-20pt}
\includegraphics[height=3in]{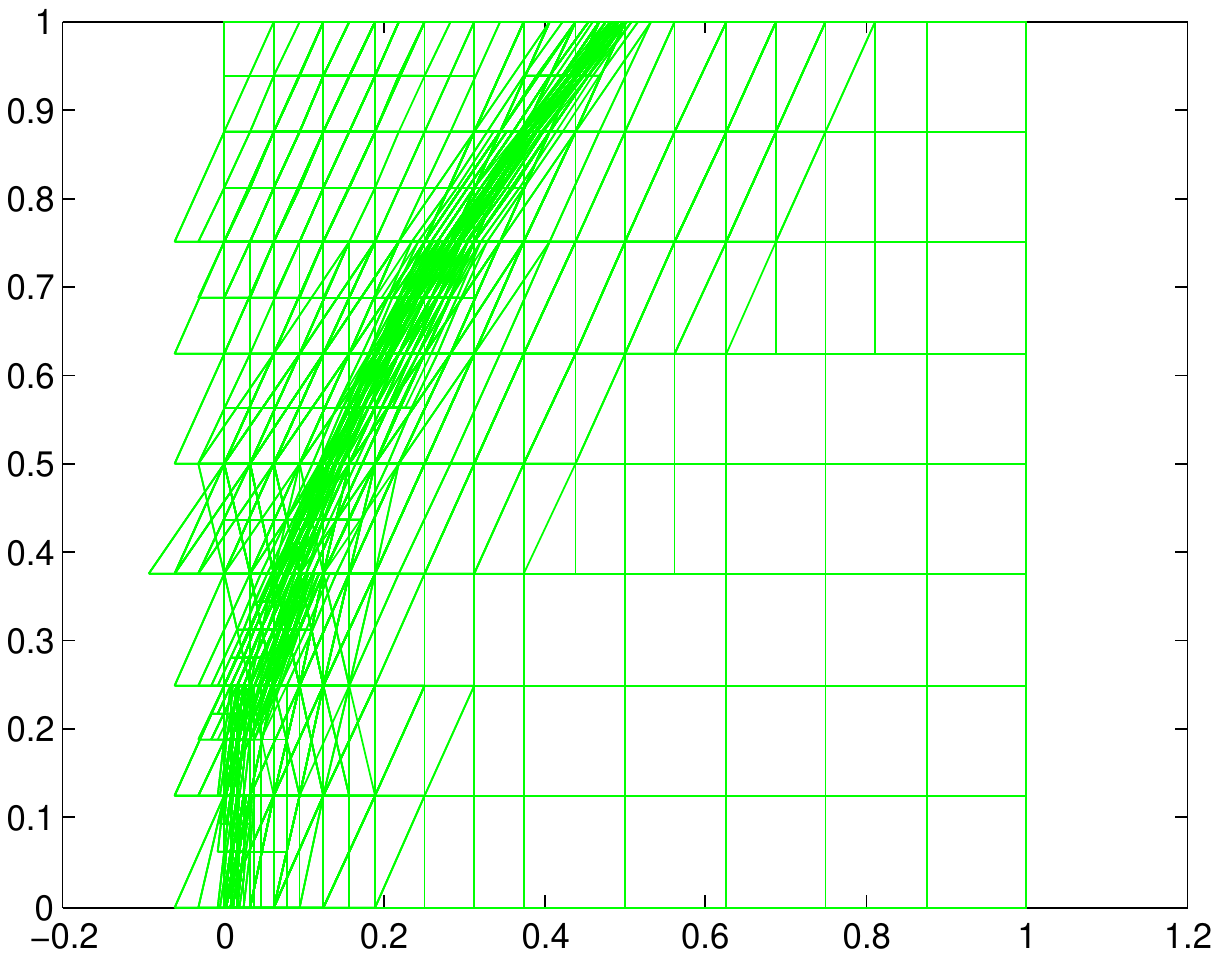}
\put(-230,60){(a)}
\put(-87,60){(b)}\\[2ex]
\vspace{-120pt}
\includegraphics[height=3in]{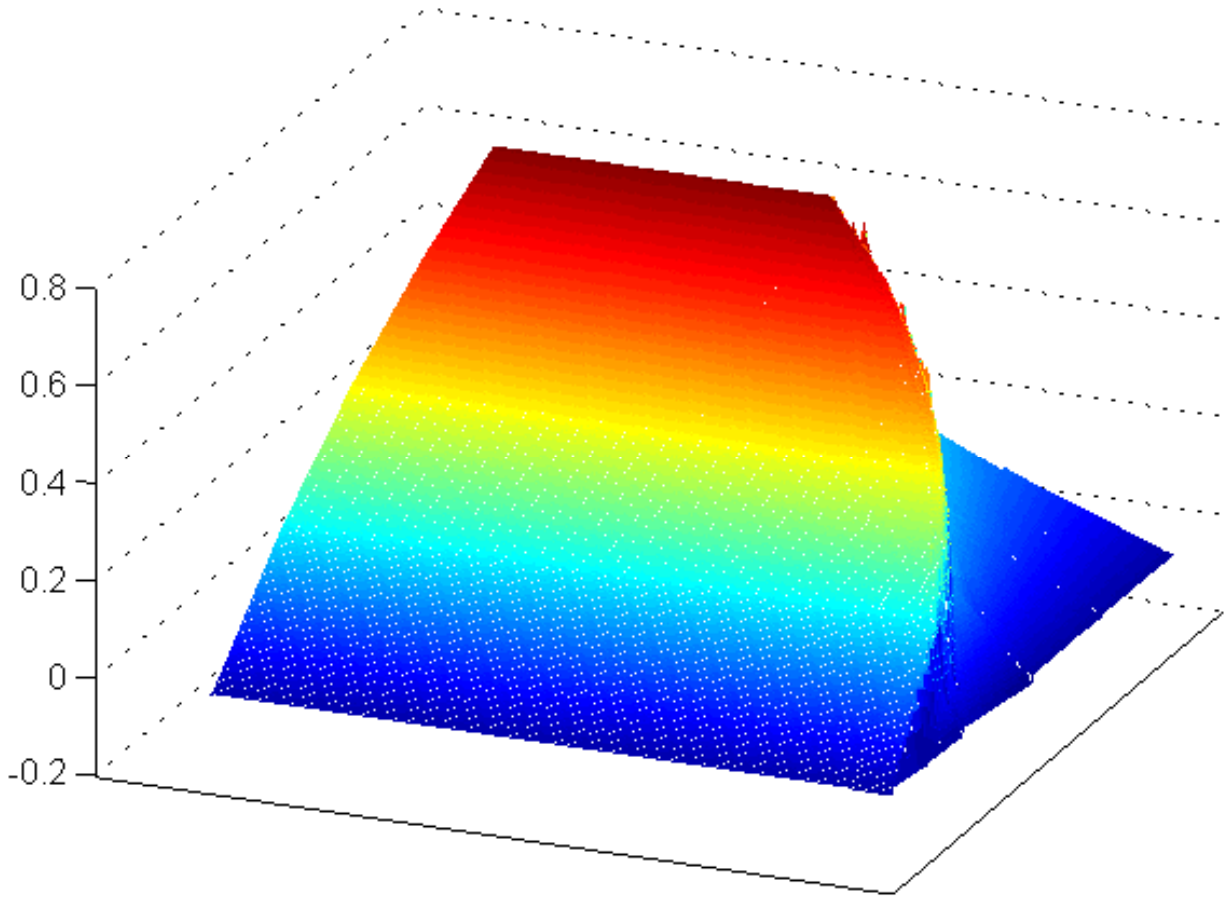}\hspace{-25pt}
\includegraphics[height=3in]{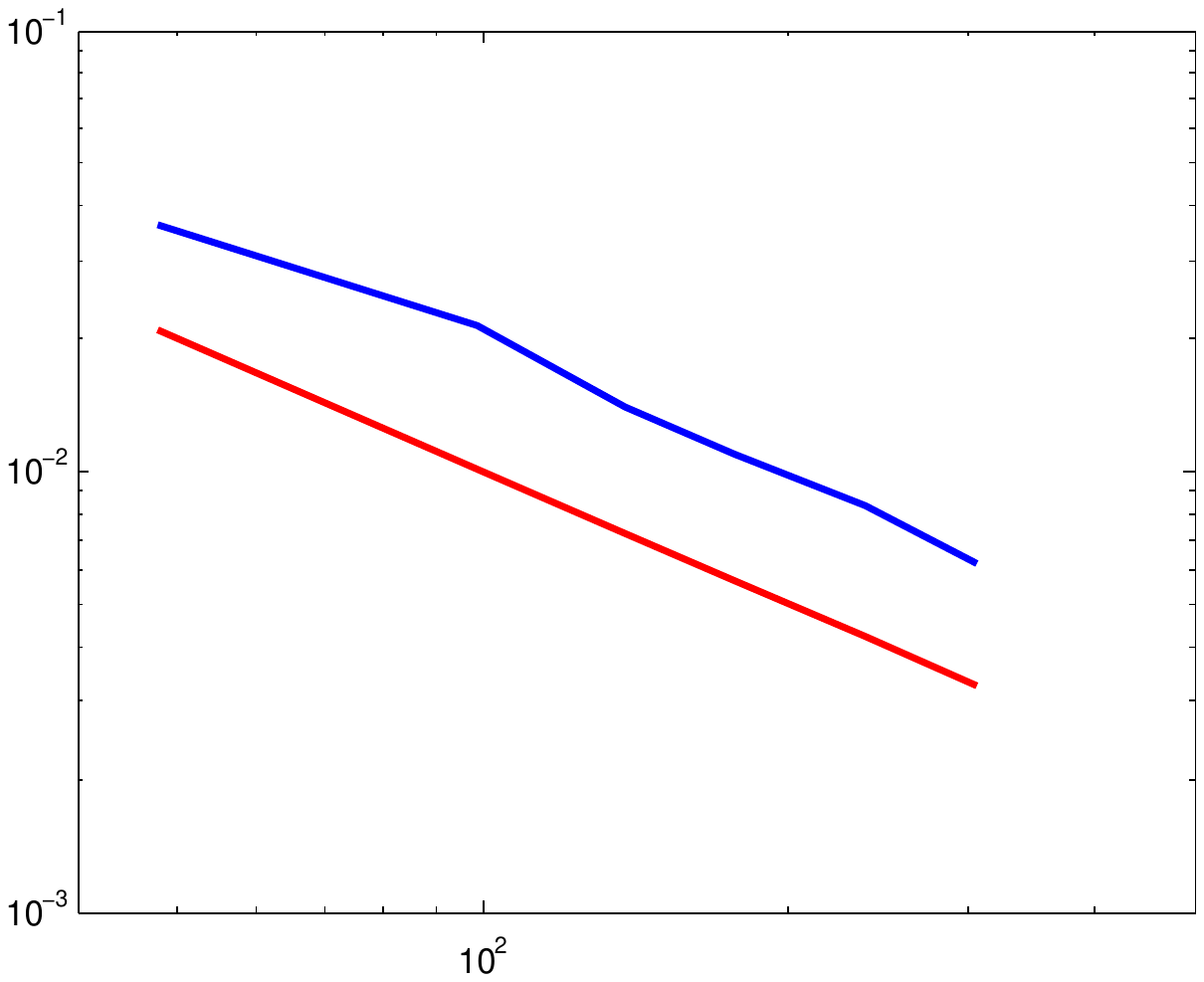}
\put(-230,60){(c)}
\put(-87,60){(d)}
\end{center}
\vspace{-50pt}
\caption{(a) Adaptive grid for the trial space $X_5$. (b) Adaptive grid for the test space $V_5$.
(c) Approximate solution (306 basis elements).
(d) $L^2(D)$ errors (vertical axis) for $N$ degrees of freedom (horizontal axis) achieved by the adaptive scheme (blue) in comparison with
the optimal  rate $N^{-1}$ (red), predicted by Theorem \ref{dahmen_huang_kutyniok_lim_schwab:ShearCart}. This is to be compared with the
rate $N^{-1/2}$ realized by adaptive {\em isotropic} refinements \cite{dahmen_huang_kutyniok_lim_schwab:DaHuScWe12}.}
\label{dahmen_huang_kutyniok_lim_schwab:shear_approx1}
\end{figure}

\vspace*{-3mm}
Figure \ref{dahmen_huang_kutyniok_lim_schwab:shear_approx1}(a), (b) show the adaptive grids
associated with the trial space $X_5$ and the test space $V_5$.
The refinement in the neighborhood of the discontinuity curve
reflects a highly anisotropic structure.
Figure \ref{dahmen_huang_kutyniok_lim_schwab:shear_approx1}(c)
illustrates the approximation given by 306 basis elements. We  emphasize that the solution is very smooth in the
vicinity of the discontinuity curve and oscillations across the jump are almost completely absent and in fact much less pronounced than
observed for isotropic discretizations. 
Figure \ref{dahmen_huang_kutyniok_lim_schwab:shear_approx1}(d)
indicates the optimal rate realized by
our scheme, see Theorem \ref{dahmen_huang_kutyniok_lim_schwab:ShearCart}.
%
%
The estimated values of the proximality parameter $\delta$, displayed in Table
\ref{dahmen_huang_kutyniok_lim_schwab:shear_delta1}, indicate the numerical stability  of the scheme.
\vspace*{-2mm}
\begin{table}
\begin{center}
\begin{tabular}{|c|c|c|}
  \hline
  $n$ & Estimated $\delta$ & $\|u_j^K - u\|_{L_2([0,1]^2)}$ \\
  \hline
   48 & 0.298138 & 0.036472 \\
  99 &  0.442948 & 0.021484 \\
  138 & 0.352767 &  0.013948 \\
  177 & 0.322156 & 0.010937 \\
  237 & 0.316545 & 0.008348 \\
  306 & 0.307965 &  0.006152 \\
  \hline
\end{tabular}
\end{center}
\caption{
Numerical estimates \eqref{eq:DeltaProx} for  the proximality constant $\delta$
and for the $L_2$ approximation error.
}
\label{dahmen_huang_kutyniok_lim_schwab:shear_delta1}
\end{table}

\vspace*{-2mm}
In the remainder of the paper 
we discuss parametric equations whose solutions are
functions of spatial variables and additional parameters.
Particular attention will here be paid to the  {\em radiative transfer problems},
where the dimension of the physical domain is $2$ or $3$.\vspace*{-6mm}
\section{Reduced Basis Methods}
\label{dahmen_huang_kutyniok_lim_schwab:RBM}
\vspace*{-3mm}
\subsection{Basic Concepts and Rate Optimality}
\vspace*{-3mm}
Model reduction is often necessary when solutions to {\em parametric families} of PDEs are frequently queried for
different parameter values e.g. in an online design or optimization process.
 The linear transport equation
 {\eqref{dahmen_huang_kutyniok_lim_schwab:rad1.1}} is a simple example of such a {\em parameter dependent}
PDE. Since   a) propagation of singularities is present and
b) the parameters determine the propagation direction $\vec{s}$ it turns out to  already pose serious difficulties for
standard model reduction techniques.


We emphasize that, rather than considering a {\em single} variational formulation for functions of spatial variables and parameters, as will be done
later in  Section \ref{dahmen_huang_kutyniok_lim_schwab:RadTrans}, we
take up the parametric nature of the problem  by considering a {\em parametric family of variational formulations}. That is, for each fixed $\vec{s}$ the problem is an ordinary linear transport problem for which we can employ the corresponding variational formulation from Section \ref{dahmen_huang_kutyniok_lim_schwab:transport},
where now the  respective spaces may depend on the parameters. In this section we summarize some of the results from \cite{dahmen_huang_kutyniok_lim_schwab:DePeWo13}
which are based in an essential way on the concepts discussed in the previous section.

In general, consider a familiy
\begin{equation}
\label{dahmen_huang_kutyniok_lim_schwab:RB1}
b_\mu(u,v) = f(v),\quad u\in X_\mu, \, v\in Y_\mu,\,\, \mu \in \mathcal{P}, \,\, b_\mu(u,v)=\sum_{k=1}^M\Theta_k(\mu)b_k(u,v)
\end{equation}
of well-posed problems, where $\mathcal{P}\subset \mathbb{R}^P$ is a compact set of parameters  {$\mu$},
and the parameter dependence
is assumed to be {\em affine} with smooth functions $\Theta_k$.
The solutions $u(\cdot;\mu)=u(\mu)$ then
become functions of the spatial variables and of the parameters $\mu\in \mathcal{P}$.

As before we can view \eqref{dahmen_huang_kutyniok_lim_schwab:RB1} as
a  {parametric} family of operator equations $B_\mu u=f$, where
$B_\mu : X_\mu \to Y_\mu'$ is again given by $(B_\mu u)(v)=b_\mu(u,v)$.
Each particular solution $u(\mu)$ is a {\em point} on the {\em solution manifold}
\begin{equation}
\label{dahmen_huang_kutyniok_lim_schwab:RBsolman}
\mathcal{M} := \{ B_\mu^{-1}f: \mu\in \mathcal{P}\}.
\end{equation}
Rather than viewing $u(\mu)$  as a point in a very high-dimensional (in fact
infinite dimensional) space, and calling  a standard solver for each evaluation in a frequent query problem,
the {\em Reduced Basis Method} (RBM) tries to exploit the fact that each $u(\mu)$ belongs
to a much smaller dimensional manifold $\mathcal{M}$.  
Assuming that all the spaces $X_\mu$ are equivalent to a reference Hilbert space $X$ with norm $\|\cdot\|_X$,  the key objective of the RBM is to construct a possibly small dimensional linear space
$X_n\subset X$ such that for a given target accuracy $\epsilon>0$
\begin{equation}
\label{dahmen_huang_kutyniok_lim_schwab:RBtarget}
\sup_{\mu \in\mathcal{P}} \inf_{w\in X_n}\|u(\mu)- w\|_X := {\rm maxdist}_X(\mathcal{M},X_n)\leq \epsilon .
\end{equation}
Once $X_n$ has been found, bounded linear functionals of the exact solution $u(\mu)$ can be approximated within accuracy
$\epsilon$ by the functional applied to an approximation from $X_n$ which, when $n$ is small, can hopefully be determined at very low cost.
The computational work in an RBM is therefore divided into an {\em offline} and an {\em online} stage. Finding $X_n$ is the core offline task
which is allowed to be computationally (very) expensive. More generally, solving problems in the ``large'' space $X$
is part of the offline stage. Of course, solving a problem in $X$ is already idealized.
In practice $X$ is replaced by a possibly very large trial space, typically a finite element space, which is referred to as the {\em truth} space
and should be chosen large  enough to guarantee the desired target accuracy, ideally certified by a posteriori bounds.

The computation of a (near-)best approximation $u_n(\mu)\in X_n$ is then to
be {\em online feasible}. More precisely, one seeks to obtain a representation
\begin{equation}
\label{dahmen_huang_kutyniok_lim_schwab:RBsep}
u_n(\mu) = \sum_{j=1}^n c_j(\mu)\phi_j,
\end{equation}
where the $\phi_j$ form a basis for $X_n$ and where for each query $\mu\in \mathcal{P}$
the expansion coefficients $c_j(\mu)$
can be computed by solving only problems of the size $n$, see e.g. \cite{dahmen_huang_kutyniok_lim_schwab:SeVeHuDeNgPa06} for principles of practical realizations.
Of course, such a concept pays off when the dimension $n=n(\epsilon)$,
needed to realize \eqref{dahmen_huang_kutyniok_lim_schwab:RBtarget}, grows very slowly when $\epsilon$ decreases.
This means that the elements of $\mathcal{M}$ have sparse representations with respect to
certain {\em problem dependent} dictionaries.

The by now most prominent strategy for constructing ``good'' spaces $X_n$ can be sketched as follows.
Evaluating for a given $X_n$ the quantity ${\rm maxdist}_X(\mathcal{M},X_n)$ is infeasible because this would require
to determine for {\em each} $\mu\in\mathcal{P}$ (or for {\em each} $\mu$ in a large training set $\mathcal{P}_h\subset \mathcal{P}$ which for simplicity  we
also denote by $\mathcal{P}$) the solution $u(\mu)$ which even for the offline stage is way too expensive.
Therefore, one chooses a {\em surrogate} $R_n(\mu)$ such that
\begin{equation}
\label{dahmen_huang_kutyniok_lim_schwab:RBsurr}
\inf_{w\in X_n}\|u(\mu)- w\|_X\leq R_n(\mu,X_n),\quad \mu\in \mathcal{P},
\end{equation}
where the evaluation of $R_n(\mu,X_n)$ is fast and an optimization
of $R_n(\mu,X_n)$ can therefore be performed in the offline stage.
This leads to the {\em greedy algorithm} in Algorithm \ref{dahmen_huang_kutyniok_lim_schwab:greedy}.
\begin{algorithm}[htb]
  \caption{greedy algorithm}
  \label{dahmen_huang_kutyniok_lim_schwab:greedy}
  \begin{algorithmic}[1]
    \Function{{\bf GA}}{}
  \State Set $X_0 := \{0\}$, $n=0$,
  \While{${\rm argmax}_{\mu\in\mathcal{P}} R(\mu,X_n) \ge \epsilon$}
    \State
    \begin{equation}
      \begin{array}{ll}
        \mu_{n+1} & := {\rm argmax}_{\mu\in\mathcal{P}} R(\mu,X_n), \\
        u_{n+1} & := u(\mu_{n+1}), \\
        X_{n+1}& := {\rm span}\,\big\{X_n,\{u(\mu_{n+1})\}\big\} = {\rm span}\,\{u_1, \dots, u_{n+1}\}
      \end{array}
      \label{dahmen_huang_kutyniok_lim_schwab:greedy1}
    \end{equation}
  \EndWhile
  \EndFunction
  \end{algorithmic}
\end{algorithm}
A natural question is to ask how the spaces $X_n$ constructed in such
a greedy fashion compare with ``best spaces''
in the sense of the {\em Kolmogorov $n$-widths}
\begin{equation}
\label{dahmen_huang_kutyniok_lim_schwab:RBkolm}
d_n(\mathcal{M})_{X} := \inf_{{\rm dim}\,W_n=n}\sup_{w\in cM}\inf_{z\in W_n}\|w-z\|_X
\;.
\end{equation}
The $n$-widths are expected to decay the faster the more regular the dependence of $u(\mu)$ is on $\mu$.
In this case an RBM has a chance to perform well.

Clearly, one always has $d_n(\mathcal{M})_{X} \leq {\rm maxdist}_X(\mathcal{M},X_n)$. Unfortunately, the best constant $C_n$
for which ${\rm maxdist}_X(\mathcal{M},X_n) \leq C_n d_n(\mathcal{M})_{X}$ is $C_n = 2^n$, see \cite{dahmen_huang_kutyniok_lim_schwab:BiCoDaDePeWo11,dahmen_huang_kutyniok_lim_schwab:BuMaPa12}.
Nevertheless, when comparing {\em rates} rather than individual values, one arrives at more positive results 
\cite{dahmen_huang_kutyniok_lim_schwab:BiCoDaDePeWo11,dahmen_huang_kutyniok_lim_schwab:DePeWo13-2}. The following 
consequence of these results asserts optimal performance of the greedy algorithm provided that the surrogate 
sandwiches the error of best approximation \cite{dahmen_huang_kutyniok_lim_schwab:DePeWo13}.

\begin{theorem}
\label{dahmen_huang_kutyniok_lim_schwab:rateoptimal}
Assume that there exists a constant $0 < c_R \leq 1$
such that for all $n$ holds
\begin{equation}
\label{dahmen_huang_kutyniok_lim_schwab:RBsurrtight}
c_R R_n(\mu,X_n) \leq \inf_{w\in X_n}\|u(\mu)- w\|_X\leq R_n(\mu,X_n),\quad \mu\in \mathcal{P}.
\end{equation}
Then, the spaces $X_n$ produced by Algorithm \ref{dahmen_huang_kutyniok_lim_schwab:greedy} satisfy
\begin{equation}
\label{dahmen_huang_kutyniok_lim_schwab:RBrates}
d_n(\mathcal{M})_x \leq C n^{-\alpha}\quad \Longrightarrow \quad {\rm maxdist}_X(\mathcal{M},X_n)\leq \bar C n^{-\alpha},
\end{equation}
where $\bar C$ depends only on $C, \alpha$, and $\kappa (R_n) := 1/c_R$,
the {\em condition} of the surrogate.
\end{theorem}

We call the RBM {\em rate-optimal} whenever \eqref{dahmen_huang_kutyniok_lim_schwab:RBrates} holds for any $\alpha >0$. Hence, finding rate-optimal
RBMs amounts to finding {\em feasible well-conditioned surrogates}.\vspace*{-6mm}

\subsection{A Double Greedy Method}
\vspace*{-4mm}
Feasible surrogates that do not require the explicit computation  of truth solutions for each $\mu\in\mathcal{P}$ need to be
based in one way or the other on {\em residuals}. When \eqref{dahmen_huang_kutyniok_lim_schwab:RB1} is a family of uniformly $X$-elliptic problems so that
$B_\mu$ are uniformly bounded isomorphisms from $X$ onto $X'$, residuals indeed lead to feasible surrogates whose condition depends
on the ratio of the continuity and coercivity constant. This follows from the mapping property of $B_\mu$, stability of the Galerkin method,
and the best approximation property of the Galerkin projection, see  \cite{dahmen_huang_kutyniok_lim_schwab:DePeWo13}.

When the problems \eqref{dahmen_huang_kutyniok_lim_schwab:RB1} are indefinite or unsymmetric
and singularly perturbed these mechanisms no longer work in this way, which explains why the conventional RBMs do not perform
well for transport dominated problems in that they are far from rate-optimal.

As shown in \cite{dahmen_huang_kutyniok_lim_schwab:DePeWo13}, a remedy is offered by the above {\em  renormation principle} providing
well-conditioned variational formulations for \eqref{dahmen_huang_kutyniok_lim_schwab:RB1}. In principle, these allow one to relate errors (in a norm of choice)
to residuals in a suitably adapted dual norm which are therefore candidates for surrogates.
The problem is that, given a trial space $X_n$, in particular a space generated in the
context of an RBM, it is not clear how to obtain a sufficiently good test space such that the corresponding Petrov-Galerkin projection
is comparable to the best approximation.
The new scheme developed in  {\cite{dahmen_huang_kutyniok_lim_schwab:DePeWo13}}
is of the following form:

\begin{itemize}
\setlength{\labelsep}{0.4cm}
\item[(I)] Initialization: take $X_1:= {\rm span}\,\{u(\mu_1)\}$, $\mu_1$ randomly chosen, $Y_1:= \{0\}$;
\item[(II)]
given a pair of spaces $X_n, \tilde V_n$, the routine \textsc{Update-inf-sup-$\delta$} enriches $\tilde V_n$ to a larger space
$V_n$ which is $\delta$-proximal for $X_n$;
\item[(III)]
extend $X_n$ to $X_{n+1}$ by a greedy step according to Algorithm \ref{dahmen_huang_kutyniok_lim_schwab:greedy}, set $\tilde V_{n+1}=V_n$, and
go to (II) as long as a given target tolerance for an a posteriori threshold is not met.
\end{itemize}

The rountine \textsc{Update-inf-sup-$\delta$}  works roughly as follows (see also \cite{dahmen_huang_kutyniok_lim_schwab:GeVe12} in the case of the Stokes system).
First, we search for a parameter $\bar{\mu} \in \mathcal{P}$ and a function $\bar{w} \in X_n$ for which the inf-sup condition   is worst, i.e.
\begin{equation}
\label{dahmen_huang_kutyniok_lim_schwab:arginfsup}
\sup_{v\in \tilde V_n}\frac{b_{\bar{\mu}}(\bar{{w}}, v)}{\|v\|_{Y_{\bar{\mu}}}\|\bar{{w}}\|_{\hat X_{\bar{\mu}}}} = \inf_{\mu\in \mathcal{P}} \left( \inf_{w\in X_n} \sup_{v\in \tilde V_n}\frac{b_\mu(w, v)}{\|v\|_{Y_\mu}\|{w}\|_{\hat X_\mu}} \right).
\end{equation}
If this worst case inf-sup constant does not exceed yet a desired uniform lower bound, $\tilde V_n$
does not contain an effective {\em supremizer}, i.e.,
a function realizing the supremum in \eqref{dahmen_huang_kutyniok_lim_schwab:arginfsup},
for $\bar{\mu}, \bar{w}$, yet.
However, since the truth space satisfies a uniform inf-sup condition, due to the
same variational formulation, there exists a good supremizer in the truth space which,
is given by the Galerkin problem
\[
\bar{v} = R_{Y_{\bar{\mu}}}^{-1}B_{\bar{\mu}} \bar{w}
=
{\rm argmax}_{v \in Y_{\bar{\mu}}}\frac{b_{\bar{\mu}}(\bar{w}, v)}{\|v\|_{Y_{\bar{\mu}}} \|\bar{w}\|_{\hat X_{\bar{\mu}}}},
\]
 providing the enrichment
$  \tilde V_n \to {\rm span} \{ \tilde V_n, R_{Y_\mu}^{-1}B_\mu  \bar{w}\}$.

The {\em interior greedy stabilization loop} {(II)} ensures that
the input pair $X_n, Y_n$ in step (III) is inf-sup stable with an inf-sup constant as close to
one as one wishes, depending on the choice of $\delta <1$.
By Theorem \ref{dahmen_huang_kutyniok_lim_schwab:stab2}, each solution $u_n(\mu)$ of
the discretized system 
for $(X_h,V) = (X_n,V_n)$ satisfies the
near-best approximation property \eqref{dahmen_huang_kutyniok_lim_schwab:BAP1}, \eqref{dahmen_huang_kutyniok_lim_schwab:BAP2}.
Hence $\|f- B_\mu u_n(\mu)\|_{Y_\mu'}$
is a well conditioned surrogate (with condition close to one).
Therefore, the assumptions of Theorem \ref{dahmen_huang_kutyniok_lim_schwab:rateoptimal} hold so that
the {\em outer greedy } step (III) yields a rate-optimal update.
In summary, under
the precise assumptions detailed in \cite{dahmen_huang_kutyniok_lim_schwab:DePeWo13},
the above {\em double greedy scheme} is {\em rate-optimal}.

Before turning to numerical examples,
a few comments on the interior greedy loop \textsc{Update-inf-sup-$\delta$}
are in order.

 (a) Finding $\bar\mu$ in \eqref{dahmen_huang_kutyniok_lim_schwab:arginfsup} requires for each $\mu$-query to perform a singular value decomposition in
 the low dimensional reduced spaces so that this is offline feasible, see \cite[Remark 4.2]{dahmen_huang_kutyniok_lim_schwab:DePeWo13}.

 (b) When the test spaces $Y_\mu$    {all agree
 with a reference Hilbert space $Y$ as sets and with with equivalent norms} it is
easy to see that the interior stabilization loop terminates after
at most $M$ steps where $M$ is the number of parametric components
in \eqref{dahmen_huang_kutyniok_lim_schwab:RB1}, see \cite[Remark 4.9]{dahmen_huang_kutyniok_lim_schwab:DePeWo13} and \cite{dahmen_huang_kutyniok_lim_schwab:GeVe12,dahmen_huang_kutyniok_lim_schwab:RoVe07}. 
If, on the other hand, the spaces $Y_\mu$ differ even as sets,
as in the case of transport equations when the transport direction is the parameter,
this is not clear beforehand.
By showing that the inf-sup condition is equivalent to a $\delta$-proximality condition one can show
under mild assumptions
though that the greedy interior loop still terminates after a number of steps which is
independent of the truth dimension, \cite[Remark 4.11]{dahmen_huang_kutyniok_lim_schwab:DePeWo13}.

(c) In this latter case the efficient evaluation of $\|f- B_\mu u(\mu)\|_{Y_\mu'}$
requires additional efforts, referred to as {\em iterative tightening},
see 
\cite[Section 5.1]{dahmen_huang_kutyniok_lim_schwab:DePeWo13}.

(d) The renormation strategy saves an expensive computation of stability constants as in conventional RBMs since, by construction, through the
choice of $\delta$, the stability constants can be driven as close to one as one wishes.

The scheme has been applied in \cite{dahmen_huang_kutyniok_lim_schwab:DePeWo13} to convection-diffusion and pure transport problems where the convection directions are parameter dependent. Hence the variational formulations are of the form \eqref{dahmen_huang_kutyniok_lim_schwab:RB1}. We briefly report some results for the transport problem, since this
is an extreme case in the following sense. The test spaces $Y_\mu$ do {\em not} agree as sets when one would like the $X_\mu$
to be equivalent for different parameters. Hence, one faces the obstructions mentioned in (b), (c) above. Moreover, for discontinuous right hand side
and discontinuous boundary conditions the dependence of the solutions on the parameters has low regularity  so that the $n$-widths do not decay
as rapidly as in the convection-diffusion case. Nevertheless, the rate-optimality still shows a relatively fast convergence for the reduced spaces $X_n$
shown below.

The first example concerns \eqref{dahmen_huang_kutyniok_lim_schwab:rad1.1} (with $\mu=\vec{s}$ ranging over a quarter circle, $D=(0,1)^2$) for $f_\circ \equiv 1$, $g\equiv 0$.
In the second example, we take $f_\circ(x_1,x_2)= 0.5$ for $x_1< x_2$, $f_\circ(x_1,x_2)= 1$ for $x_1\geq x_2$. 

\begin{table}
  \centering
  \begin{tabular}{|c|c|c|c|c|c|c|}
      \hline
      \multicolumn{2}{|c|}{dimension} & & maximal & \multicolumn{2}{|c|}{maximal error between} & surr / \\
      trial & test & $\delta$ & surr &  rb truth &  rb $L_2$ & err \\
      \hline
      \hline
      4 & 11 & 3.95e-01 & 8.44e-03 & 2.45e-02 & 2.45e-02 & 3.45e-01 \\
      10 & 33 & 4.32e-01 & 3.37e-03 & 5.74e-03 & 5.74e-03 & 5.87e-01 \\
      16 & 57 & 4.32e-01 & 1.50e-03 & 2.56e-03 & 2.56e-03 & 5.84e-01 \\
      20 & 74 & 4.16e-01 & 1.21e-03 & 2.10e-03 & 2.10e-03 & 5.77e-01 \\
      24 & 91 & 4.05e-01 & 7.27e-04 & 1.58e-03 & 1.58e-03 & 4.61e-01 \\
      \hline
  \end{tabular}
  \caption{Numerical results for Example 1, maximal $L_2$ truth error 0.000109832.}
  \label{dahmen_huang_kutyniok_lim_schwab:transport-1}
\end{table}
\vspace*{-7mm}

\begin{table}
  \centering
  \begin{tabular}{|c|c|c|c|c|c|c|}
      \hline
      \multicolumn{2}{|c|}{dimension} & & maximal & \multicolumn{2}{|c|}{maximal error between} & surr / \\
      trial & test & $\delta$ & surr & rb truth & rb $L_2$ & err \\
      \hline
      \multicolumn{7}{|c|}{first reduced basis creation} \\
      \hline
      20 & 81 & 3.73e-01 & 2.71e-02 & 5.46e-02 &  5.62e-02 & 4.82e-01 \\
      \hline
      \multicolumn{7}{|c|}{second reduced basis creation} \\
      \hline
      10 & 87 & 3.51e-01 & 6.45e-02 & 7.40e-02 &  7.53e-02 & 8.57e-01 \\
      \hline
  \end{tabular}
  \caption{Numerical results for Example 2 after a single cycle of iterative tightening.
  Maximal $L_2$ truth error 0.0154814}
\end{table}

\vspace*{-1cm}

\begin{figure}[htb]
  \begin{center}
    \includegraphics[width=0.45\textwidth]{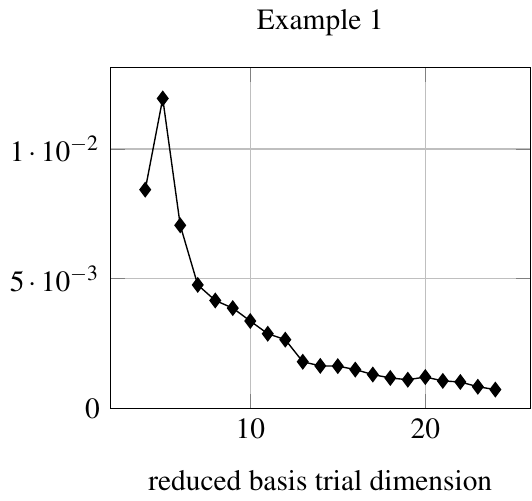}
    \includegraphics[width=0.45\textwidth]{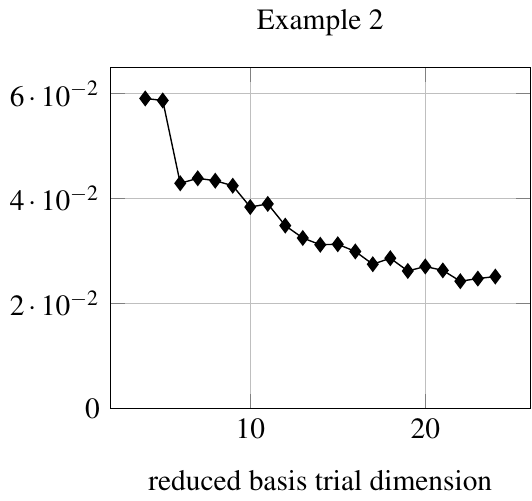}
  \end{center}
  \vspace*{-3mm}
  \caption{Surrogates of the reduced basis approximation for Examples 1 and 2.}
\end{figure}

\vspace*{-1cm}

\section{Sparse Tensor Approximation for Radiative Transfer}
\label{dahmen_huang_kutyniok_lim_schwab:RadTrans}
\vspace*{-3mm}
 We now extend the parametric transport problem \eqref{dahmen_huang_kutyniok_lim_schwab:rad1.1} to  the {\em radiative transport problem} (RTP) (see, eg., \cite{dahmen_huang_kutyniok_lim_schwab:Mo03})
which consists in finding the \emph{radiative intensity}
$u : D \times \mathcal{S} \rightarrow \mathbb{R}$,
defined on the Cartesian product of a bounded physical domain
$D\subset \mathbb{R}^d$, where $d=2,3$,
and the unit $d_\mathbb{S}$-sphere as the parameter domain:
$\mathcal{P} = \mathcal{S}$ with $d_\mathbb{S} = 1,2$.
Given an  \emph{absorption coefficient} $\kappa \geq 0$,  a \emph{scattering coefficient} $\sigma \geq 0$, and
a \emph{scattering kernel} or \emph{scattering phase function} $\Phi > 0$, which
is normalized to
$\int_{\mathcal{S}} \Phi(\vec{s}, \vec{s}') {\rm d}{\vec{s}'} = 1$
for each direction $\vec{s}$, one defines
the transport operator $\mathrm{T} u := (\vec{s}\cdot{\boldmath \nabla}_x + \kappa) u$,
and the scattering operator
$\mathrm{Q} u := \sigma \mathrm{Q_1}\!
u = \sigma (u - \int_{\mathcal{S}} \Phi(\vec{s}, \vec{s}')u(\vec{x},\vec{s}')
{\rm d}{\vec{s}'})$.
The radiative intensity is then given by
\begin{equation}
\label{dahmen_huang_kutyniok_lim_schwab:RadTransEq}
(\mathrm{T} + \mathrm{Q}) u = f, \;\;
u|_{\partial \Omega_- } = g,
\end{equation}
where $ f := \kappa I_b$,
$\partial \Omega_-
:=
\{(\vec{x}, \vec{s}) \in \partial D \times \mathcal{S} : \vec{s} \cdot \vec{n}(\vec{x}) < 0\}$,
and $g$ denote the source term,  the inflow-boundary,
and the inflow-boundary values, respectively.
As before, $\Gamma_-(\vec{s}) := \{\vec{x} \in \partial D : \vec{s} \cdot \vec{n}(\vec{x}) < 0\}$
stands for the physical inflow-boundary.
%
%
%

The partial differential equation~\eqref{dahmen_huang_kutyniok_lim_schwab:RadTransEq}
is known as \emph{stationary monochromatic radiative transfer equation} (RTE)
with scattering, and can be viewed as (nonlocal) extension of the parametric
transport problem \eqref{dahmen_huang_kutyniok_lim_schwab:rad1.1},
where the major difference  {to \eqref{dahmen_huang_kutyniok_lim_schwab:rad1.1}
is the scattering operator $Q$.}
Sources with support contained in $D$ are modeled by the \emph{blackbody intensity}
$I_b \geq 0$, radiation from sources outside of the domain or
from its enclosings is prescribed by the \emph{boundary data} $g \geq 0$.
The vector $\vec{n}(\vec{x})$ denotes the \emph{outer unit normal}
on the boundary $\partial D$ of the physical domain.
%
%

Deterministic numerical methods for the RTP which are commonly used
in engineering comprise the
{\em discrete ordinates ($S_N$-) method} and the
{\em spherical harmonics ($P_N$-) method}.

In the {\em discrete ordinate method} (DOM),
the angular domain is collocated by a finite number of fixed
propagation directions in the angular parameter space;
in this respect, the DOM   resembles the
greedy collocation in the parameter domain:
each of the directions Eq.~\eqref{dahmen_huang_kutyniok_lim_schwab:RadTransEq}
results in a spatial PDE which is solved (possibly in parallel)
by standard finite differences, finite elements, or finite volume methods.

In the {\em spherical harmonics method} (SHM),
a spectral expansion with spatially variable coefficients
is inserted as ansatz into the variational principle
Eq.~\eqref{dahmen_huang_kutyniok_lim_schwab:RadTransEq}.
By orthogonality relations, a coupled system of PDEs (whose type can change from
hyperbolic to elliptic in the so-called diffuse radiation approximation)
for the spatial coefficients is obtained,
which is again solved by finite differences or finite elements.

The common deterministic methods $S_N$- and $P_N$-approximation exhibit the so-called ``curse of dimensionality'':
the error with respect to the total numbers of degrees of freedom (DoF) $M_D$ and $M_\mathcal{S}$
on the physical domain $D$ and the parameter domain
$\mathcal{S}$ scales with the dimension
$d$ and $d_\mathbb{S}$ 
as $O(M_D^{-s/d} + M_\mathcal{S}^{-t/d_\mathbb{S}})$ with
positive constants $s$ and $t$.

The so called {\em sparse grid approximation method} alleviates this curse of
dimensionality for elliptic PDEs on cartesian product domains,
see \cite{dahmen_huang_kutyniok_lim_schwab:BuGr04} and the references therein.
\cite{dahmen_huang_kutyniok_lim_schwab:WiHiSc08} has developed a sparse tensor method to overcome
the curse of dimensionality for radiative transfer with a
wavelet (isotropic) discretization of the angular domain.
Under certain regularity assumptions on the absorption coefficient $\kappa$ and the blackbody intensity $I_b$, their method achieves the typical benefits of sparse tensorization: a log-linear complexity in the number of degrees of freedom of a component domain with an essentially (up to a logarithmic factor) undeteriorated rate of convergence.
However, scattering had not been addressed in that work.

In order to include scattering and to show that the concepts of sparse tensorization
can also be applied to common solution methods,
sparse tensor versions of the spherical harmonics approximation were developed
extending the ``direct sparse'' approach by \cite{dahmen_huang_kutyniok_lim_schwab:WiHiSc08}.
The presently developed version also accounts for scattering \cite{dahmen_huang_kutyniok_lim_schwab:GrSc11a}.
For this sparse spherical harmonics method, we proved that
the benefits of sparse tensorization can indeed be harnessed.
%

As a second method a {\em sparse tensor product version of the DOM}
based on the {\em sparse grid combination technique} was realized and analyzed
in \cite{dahmen_huang_kutyniok_lim_schwab:Gr13,dahmen_huang_kutyniok_lim_schwab:GrSc11b}.
Solutions to discretizations of varying discretization levels,
for a number of collocated transport problems, and {\em with scattering discretized
by combined Galerkin plus quadrature approximation in the transport collocation directions}
are combined in this method to form a sparse tensor solution that
we proved in \cite{dahmen_huang_kutyniok_lim_schwab:Gr13,dahmen_huang_kutyniok_lim_schwab:GrSc11b}
breaks the curse of dimensionality as described above.
These benefits hold as long as the exact solution of the RTE is sufficiently regular.
An overview follows.\vspace*{-7mm}
\subsection{Sparse discrete ordinates method (Sparse DOM)}
\label{dahmen_huang_kutyniok_lim_schwab:DOM}
\vspace*{-3mm}
We adopt a formulation where the inflow boundary conditions
are enforced in a weak sense.
To this end, we define the boundary form
(see, eg., \cite{dahmen_huang_kutyniok_lim_schwab:Gr13})
\begin{equation}
\label{dahmen_huang_kutyniok_lim_schwab:DefBilFormBoundary}
\partial b (u, v)
:=
{(v,\vec{s}\cdot \vec{n} u)}_{L^2(\partial \Omega_-)}
=
\int_{\mathcal{S}}\int_{\Gamma_-(\vec{s})} \vec{s}\cdot\vec{n}uv{\rm d}{\vec{x}}{\rm d}{\vec{s}}
\;.
\end{equation}
For $v : D \times \mathcal{S} \rightarrow \mathbb{R}$, the norms
\begin{equation*}
\|v\|_-^2 := -\partial b(v,v), \quad {\|v\|_1}^2 := \|v\|^2 + \|\vec{s}\cdot{\boldmath \nabla}_x v\|^2 + \|\mathrm{Q_1}\! v\|^2 + \|v\|_-^2
\end{equation*}
define the Hilbert space
$\mathcal{V}_1 := \{v \in L^2(D \times \mathcal{S}) : \|v\|_1 < \infty\}$
\;.
The SUPG-stabilized Galerkin variational formulation reads:
find $u \in \mathcal{V}_1$ such that
\begin{equation}
\label{dahmen_huang_kutyniok_lim_schwab:VarFormDOM}
 ({\mathrm{R} v},{(\mathrm{T} + \mathrm{Q}) u})_{L^2({D \times \mathcal{S}})} - 2\partial b(u,v)
=
({\mathrm{R} v},{f})_{L^2(D \times \mathcal{S})} - 2\partial b(g,v) \quad \forall v \in \mathcal{V}_1
\end{equation}
with SUPG stabilization $\mathrm{R} v := v + \delta \vec{s}\cdot{\boldmath \nabla}_x v$,
where $\delta \approx 2^{-L}$.

For the discretization of \eqref{dahmen_huang_kutyniok_lim_schwab:VarFormDOM},
we replace $\mathcal{V}_1$ by $V^{L,N} = V_{D}^L \otimes V_{\mathcal{S}}^N$.
In the physical domain, standard $P_1$-FEM with a one-scale basis
on a uniform mesh of width $h_{\max} \lesssim 2^{-L}$ is used,
in the angular domain, piecewise constants on a quasiuniform
mesh of width $h_{\max} \lesssim N^{-1}$.
Fully discrete problems are obtain with
a one-point quadrature in the angular domain.
The resulting Galerkin formulation \eqref{dahmen_huang_kutyniok_lim_schwab:VarFormDOM} can be shown to result
in the same linear system of equations
as the standard collocation discretization \cite[Sec.~5.2]{dahmen_huang_kutyniok_lim_schwab:Gr13}.
The solution is constructed with
the sparse grid combination technique (cp. \cite{dahmen_huang_kutyniok_lim_schwab:BuGr04}:
$$
\hat{u}_{L,N}
= \sum_{\ell_D=0}^{L}\left(u_{\ell_D, \ell_\mathcal{S}^{\max}(\ell_D)}- u_{\ell_D,\ell_\mathcal{S}^{\max}(\ell_D+1)}\right),
$$
where $u_{\ell_D, \ell_\mathcal{S}} \in V^{\ell_D, \ell_\mathcal{S}}$
denotes the solution to a full tensor subproblem of physical resolution level
$\ell_D$ and angular resolution level $\ell_\mathcal{S}$.
The maximum angular index
$\ell_\mathcal{S}^{\max} = 2^{\lfloor \log_2(N+1)\rfloor /L (L-l_D)}$
ensures that the angular resolution decreases
when the physical resolution increases and vice versa.

While the full tensor solution $u_{L, N}$ requires
$O(2^{dL} N^{d_\mathbb{S}})$ degrees of freedom,
the sparse solution involves asymptotically at most
$O((L + \log N) (2^{dL} + N^{d_\mathbb{S}}))$ degrees
of freedom \cite[Lemma 5.6]{dahmen_huang_kutyniok_lim_schwab:Gr13}.
At the same time,
$$
\|u-u_{L,N}\|_{1}
\lesssim
2^{-L} \|u\|_{H^{2,0}(D\times \mathcal{S})} + N^{-1} \|u\|_{H^{1,1}(D\times \mathcal{S})},
$$
while for solutions in
$H^{2,1}(D\times \mathcal{S}) \subset ( H^{2,0}(D\times \mathcal{S}) \cap H^{1,1}(D\times \mathcal{S}) )$
$$
\|u-\hat{u}_{L,N}\|_{1} \lesssim L\max\{2^{-L}, N^{-1}\} \|u\|_{H^{2,1}(D\times \mathcal{S})}
\;.
$$
\vspace*{-9mm}
\subsection{Numerical experiment}
\label{dahmen_huang_kutyniok_lim_schwab:NumExRT}
\vspace*{-4mm}
To evaluate the solution numerically we monitor
the incident radiation
$G(\vec{x}) = \int_\mathcal{S} u(\vec{x}, \vec{s}) \mbox{d}\vec{s}$
and its relative error
$err(G_{L,N})_X = \| G - G_{L,N} \|_X / \|G\|_X$,
$X = L^2(D), H^1(D)$.

The setting for the experiment is $D = [0,1]^d$, $\mathcal{S} = \mathcal{S}^d_\mathbb{S}$.
We solve the RTP \eqref{dahmen_huang_kutyniok_lim_schwab:RadTransEq}
with isotropic scattering $\Phi(\vec{s}, \vec{s}') = 1/|\mathcal{S}|$ and with
zero inflow boundary conditions $g = 0$.
A blackbody radiation $I_b(\vec{x}, \vec{s})$ corresponding to the exact solution
 \[
  u(\vec{x}, \vec{s}) = \frac{3}{16 \pi}(1 + (\vec{s}\cdot \vec{s}')^2) \prod_{i=1}^3 (-4 x_i(x_i-1)),
\]
with fixed $\vec{s}' = (1/\sqrt{3}, 1/\sqrt{3}, 1/\sqrt{3})^\top$
is inserted in the right hand side functional in~\eqref{dahmen_huang_kutyniok_lim_schwab:VarFormDOM}.
The absorption coefficient is set to $\kappa = 1$,
the scattering coefficient to $\sigma = 0.5$.

This $3+2$-dimensional problem was solved with a
parallel C++ solver designed for the sparse tensor
solution of large-scale radiative transfer problems.
\begin{figure}
 \centering
\includegraphics[width=0.6\textwidth]{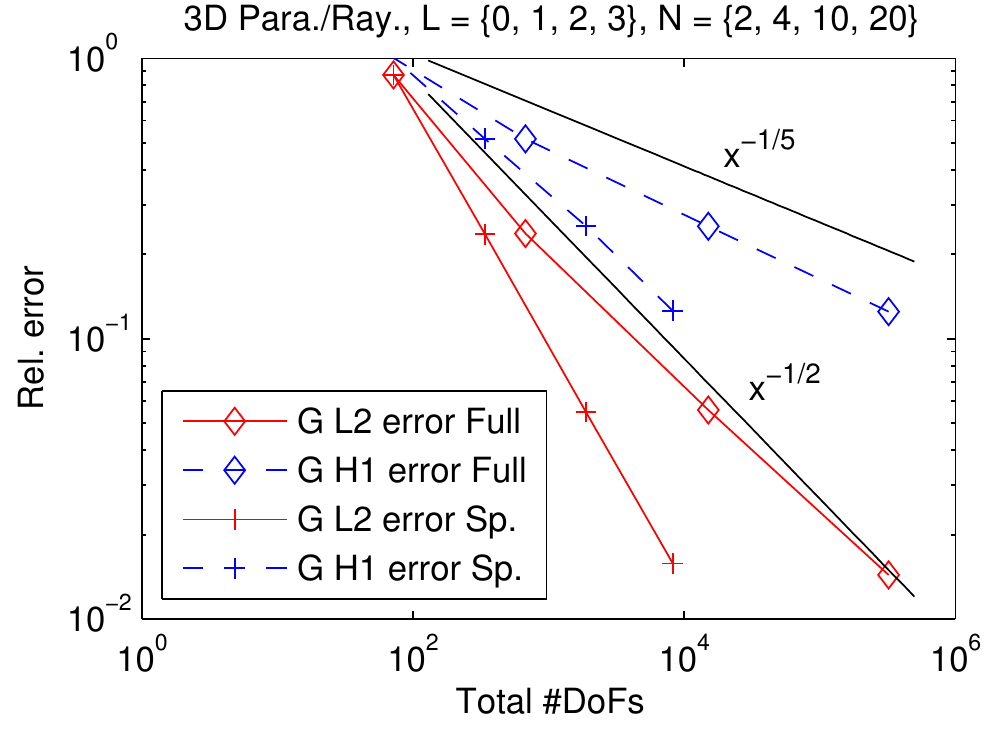}
 \caption{Convergence in incident radiation with full and sparse DOM.
Resolution for reference solution
was $L_{\mbox{ref}} = 4$.
Reference slopes provided as visual aids only.
}
\label{dahmen_huang_kutyniok_lim_schwab:ExpParaRayTDError}
\end{figure}
Fig.~\ref{dahmen_huang_kutyniok_lim_schwab:ExpParaRayTDError} shows the superior efficiency of
the sparse approach with respect to number of degrees of freedom vs.\ achieved error.
The convergence rates indicate that the curse of dimensionality
is mitigated by the sparse DOM.
Further gains are expected once the present, nonadaptive sparse DOM is
replaced by the greedy versions outlined in Section \ref{dahmen_huang_kutyniok_lim_schwab:AbsPG}.


\end{document}